\documentclass[11pt]{article}
\usepackage{amsthm}
\usepackage{amssymb}
\usepackage{amsmath}

\numberwithin{equation}{section}
\begin{document}
\title{Partial Rigidity of CR Embeddings of Real Hypersurfaces into Hyperquadrics with Small Signature Difference}
\author{PETER EBENFELT, RAVI SHROFF}
\date{}
\maketitle

\begin{abstract}
We study the rigidity of holomorphic mappings from a neighborhood of a Levi-nondegenerate CR hypersurface $M$ with signature $l$ into a hyperquadric $Q_{l'}^{N} \subseteq \mathbb{CP}^{N+1}$ of larger dimension and signature.  We show that if the CR complexity of $M$ is not too large then the image of $M$ under any such mapping is contained in a complex plane with dimension independent of $N$.  This result follows from two theorems, the first demonstrating that for sufficiently degenerate mappings, the image of $M$ is contained in a plane, and the second relating the degeneracy of mappings into different quadrics. 
\end{abstract}

\section{Introduction}

The phenomenon of rigidity of mappings between hypersurfaces embedded in complex space has been studied for many years, beginning with Poincar\'e in the early 20th century.  Initially, much work was devoted to mappings between spheres, and rigidity results were proved given restrictions on the codimension.  However, there are differences between the study of mappings between strictly pseudoconvex hypersurfaces and those assumed to be Levi nondegenerate with positive signature.  There are also differences when the source manifold is not itself a quadric but rather assumed to be embeddable into a quadric (the so-called \emph{CR complexity} is the difference between the CR dimensions of the target quadric and the source, where the CR dimension of the target is assumed to be minimal).

By the Lewy extension theorem, in the case of positive signature we need only consider restrictions of holomorphic maps.  In the case of low CR complexity but no signature difference, the main result in \cite{BEH1} says that embeddings are unique up to automorphisms of the target quadric (so-called ``super-rigidity").  In the case of zero CR complexity but positive signature difference, the main result in \cite{BEH2} says that the image of an embedding must be contained in a complex plane with dimension related to the signature difference (so-called ``partial rigidity").  In this paper we consider both low CR complexity and positive signature difference and prove a partial rigidity result.  Our proofs make use of the theory of pseudo-hermitian and pseudo-conformal geometry, particularly the work of Chern and Moser, and subsequent work of Webster.  We use recent derivations in \cite{EHZ} and \cite{BEH1} extensively.
 
Let $M \subseteq \mathbb{C}^{n+1}$ be a smooth connected Levi-nondegenerate hypersurface and $\mathcal{L}$ a representative of the Levi form of $M$.  If $M$ is connected then let $l \leq n/2$ denote the minimum of the number of positive and negative eigenvalues of $\mathcal{L}$ at any point.  This integer is constant over $M$ and will be referred to as the \emph{signature} of $M$.  We let $Q_{l}^{N} \subseteq \mathbb{CP}^{N+1}$ denote the standard hyperquadric given in homogeneous coordinates $[z_{0}:z_{1}:\ldots:z_{N+1}]$ by
\begin{equation*}
-\sum_{j=0}^{l}|z_{j}|^{2} + \sum_{k=l+1}^{N+1}|z_{k}|^{2} = 0.
\end{equation*}
Notice that $Q_{l}^{N}$ is a connected Levi-nondegenerate CR hypersurface of CR dimension $N$ and signature $l$.

We first generalize Theorem 2.2 in \cite{EHZ} which deals with degenerate smooth CR-immersions of a CR-hypersurface into a sphere.  A CR-immersion is degenerate if the span of the second fundamental form and its covariant derivatives fail to be the whole normal space of the embedding.  This and other important notions will be made precise in sections 3 and 4 of this paper.  We allow degenerate immersions into hyperquadrics where the signature of the Levi form of the target quadric is allowed to be strictly greater than that of the source manifold.  Our first result is the following

\newtheorem{Theorem 1}{Theorem}[section]
\begin{Theorem 1}
Let $M \subset \mathbb{C}^{n+1}$ be a smooth connected Levi-nondegenerate hypersurface of signature $l \le n/2$ and $f:M \rightarrow Q^{N}_{l'}$ a smooth CR mapping that is CR-transversal to $Q^{N}_{l'}$ at $f(p)$ for $p \in M$.  Assume that $f$ is constantly $(k,s)$-degenerate near $p$ for some $k$ and $s$.  If $N - n - s < n,$ then there is an open neighborhood $V$ of $p$ in $M$ such that $f(V)$ is contained in the intersection of $Q^{N}_{l'}$ with a complex plane $P \subset \mathbb{C}^{N+1}$ of codimension $s$.
\end{Theorem 1}

The idea of the proof of Theorem 1.1 goes back to the arguments in section 9 of \cite{EHZ}.

Our main result, Theorem 1.2, now follows from Theorem 1.1 and Theorem 6.1, which relates the dimensions of spaces of covariant derivatives of the second fundamental form for different embeddings.

\newtheorem{Theorem 3}[Theorem 1]{Theorem}
\begin{Theorem 3}
Let $M\subseteq \mathbb{C}^{n+1}$ be a smooth connected Levi-nondegenerate CR hypersurface with signature $l\leq n/2$.  Suppose there is an open connected neighborhood $U$ of $M$ in $\mathbb{C}^{n+1}$ and a holomorphic mapping $f_{0}:U\rightarrow \mathbb{CP}^{N_{0}+1}$ with $f_{0}(M) \subseteq Q_{l}^{N_{0}}$ and $f_{0}$ CR transversal to $Q_{l}^{N_{0}}$ along $M$.  Let $f:U \rightarrow \mathbb{CP}^{N+1}$ be a holomorphic mapping with $f(M) \subseteq Q_{l'}^{N}$ and $f$ CR transversal to $Q_{l'}^{N}$ along $M$.  Then the following hold
\begin{description}
\item[(a)] If $l = n/2$ or $f$ is side preserving then $l' \ge l$ and $N - l' \ge n - l$.  If either 
\begin{description}
\item[(i)] $(N_{0} - n) + (l' - l) < l$ or
\item[(ii)] $(N_{0} - n) + min(l' - l, (N-l')-(n-l)) < n$ and $(N - l') - (n - l) < l$, 
\end{description}
then $f(M) \subseteq Q_{l'}^{N} \cap P$, where $P \subseteq \mathbb{C}^{N+1}$ is a complex plane of dimension $(N_{0} + 1) + min(l' - l,(N-l')-(n-l))$. 

\item[(b)] If $f$ is side reversing then $N - l' \ge l$ and $l' \ge n - l$.  If $l' < n$ and $(N_{0} - n) + (l' + l -n) < n$ then $f(M) \subseteq Q_{l'}^{N} \cap P$, where $P \subseteq \mathbb{C}^{N+1}$ is a complex plane of dimension $(N_{0} + 1) + (l' + l - n)$.

\end{description}
\end{Theorem 3}

\theoremstyle{definition}
\newtheorem{remark 1}[Theorem 1]{Remark}
\begin{remark 1}
We observe that if $l' = N/2$ then the inequalities $(N- l') - (n - l) < l$ and $l' < n$ are equivalent and the  conclusions of parts \textbf{(a)} and \textbf{(b)} of Theorem 1.2 coincide.  We also observe that if $f$ is side preserving, either assumption \textbf{(i)} or \textbf{(ii)} could apply.  For instance, if $n=5$, $l = 1$, $N = 7$, $l' = 3$, and $N_{0} = 6$, then assumption \textbf{(i)} does not hold, but assumption \textbf{(ii)} does.  However, if $N$ is sufficiently large \textbf{(i)} may hold but not \textbf{(ii)}.
\end{remark 1}

\theoremstyle{definition}
\newtheorem{remark 2}[Theorem 1]{Remark}
\begin{remark 2}
Note that this partially generalizes Theorem 1.1 from \cite{BEH1} by allowing a positive signature difference between the source manifold and target hyperquadric $Q_{l'}^{N}$.  We observe that the first conclusion of part  \textbf{(b)} of the theorem implies that side reversing maps cannot exist where there is no signature difference between the source and target.  The above result also partially generalizes Theorem 1.1 from \cite{BEH2} by allowing the source manifold to be embeddable into a hyperquadric rather than be a hyperquadric itself.  The proof is given in section 6 following the statement of Theorem 6.1.
\end{remark 2}

\section{Two Important Lemmas}

We now state two key lemmas that are ingredients in the proofs of subsequent theorems.  The first lemma was proved in \cite{Hu99}.  We use the Einstein summation convention in the rest of this paper except where otherwise indicated.

\theoremstyle{plain}
\newtheorem{Huang's Lemma}{Lemma}[section]
\begin{Huang's Lemma}
Let $g_{1},\ldots,g_{k},f_{1},\ldots,f_{k}$ be holomorphic functions in $z\in\mathbb{C}^{n}$ near $0$.  Assume $g_{j}(0)=f_{j}(0)=0$ for all $j$.  Let $A(z,\bar{z})$ be real-analytic near the origin such that
\[\sum_{j=1}^{k}g_{j}(z)\overline{f_{j}(z)} = A(z,\bar{z}) (h_{a\bar{b}}z^{a}z^{\bar{b}}) \]
where $H = (h_{a\bar{b}})$ is a constant invertible matrix.  If $k < n$, then  $A(z,\bar{z}) \equiv 0.$
\end{Huang's Lemma}

Although the statement of Lemma 2.1 in \cite{Hu99} is for $H = I$, the proof for arbitrary constant invertible $H$ is identical.  We shall also need the following.

\newtheorem{Lemma 2.2}[Huang's Lemma]{Lemma}
\begin{Lemma 2.2}
Let $k,l,n$ be nonnegative integers with $k < l \leq n/2$.  Assume that $g_{1},\ldots,g_{k}, f_{1}\ldots f_{m}$ are germs at $0 \in \mathbb{C}^{n}$ of holomorphic functions and $A(z,\bar{z})$ be real-analytic near the origin such that
\[ -\sum_{i=1}^{k} |g_{i}(z)|^{2} + \sum_{j=1}^{m}|f_{j}(z)|^{2} = A(z,\bar{z})\bigg(-\sum_{i=1}^{l}|z_{i}|^{2} + \sum_{j=l+1}^{n}|z_{j}|^{2}\bigg).\]
Then $A(z,\bar{z}) \equiv 0$.
\end{Lemma 2.2}
The proof of Lemma 2.2 can be found in Lemma 4.1 of \cite{BH05} (with $l' = l$ and after an application of Lemma 2.1 of \cite{BH05}).

\section{Preliminaries}

We will use the notation of \cite{EHZ}.  Let $M$ be a Levi-nondegenerate CR-manifold of dimension $2n+1$, with rank $n$ CR bundle $\mathcal{V}$ and signature $l \leq n/2$.  Near a point $p_{0}$, we let $\theta$ be a contact form and $T$ its characteristic (or Reeb) vector field, so T is the unique real vector field satisfying $T \lrcorner d\theta = 0$ and $\langle \theta, T \rangle = 1$.  We complete $\theta$ to an admissible coframe $(\theta, \theta^{1},\ldots,\theta^{n})$ for the bundle $T'M$ of $(1,0)$-cotangent vectors (the cotangent vectors that annihilate $\mathcal{V}$.  The coframe is called admissible if $\langle \theta^{\alpha}, T \rangle = 0$, for $\alpha = 1, \ldots, n$.
We choose a frame $L_{1},\ldots,L_{n}$ for the bundle $\mathcal{\bar{V}}$ such that $(T, L_{1},\ldots,L_{n},L_{\bar{1}},\ldots,L_{\bar{n}})$ is a frame for $\mathbb{C}TM$ dual to the coframe $(\theta, \theta^{1},\ldots,\theta^{n},\theta^{\bar{1}},\ldots,\theta^{\bar{n}})$.  We use the notation that $L_{\bar{\alpha}} = \bar{L}_{\alpha}$, etc.  Relative to this frame, let $(g_{\alpha\bar{\beta}})$ denote the matrix of the Levi form.   Although we generally won't explicitly use this fact, we may assume $g_{\alpha\bar{\beta}}$ is constant and diagonal, with diagonal elements $\pm 1$ corresponding to the signature.

We denote by $\nabla$ the Tanaka-Webster connection, given relative to the chosen frame and coframe by
\begin{equation*}
\nabla L_{\alpha} := \omega_{\alpha}^{\:\:\beta}\otimes L_{\beta}.
\end{equation*}
The connection 1-forms $\omega_{\alpha}^{\:\:\beta}$ are completely determined by the conditions
\begin{align}
d\theta^{\beta} &= \theta^{\alpha}\wedge\omega_{\alpha}^{\:\:\beta} \qquad \textrm{mod}\theta\wedge\theta^{\bar{\alpha}}, \nonumber \\
dg_{\alpha\bar{\beta}} &= \omega_{\alpha\bar{\beta}} + \omega_{\bar{\beta}\alpha}.
\end{align}
Note that we use the Levi form to lower and raise indices as usual.  We may rewrite the first condition in (3.1) as 
\begin{equation}
d\theta^{\beta} = \theta^{\alpha}\wedge\omega_{\alpha}^{\:\:\:\beta} + \theta\wedge\tau^{\beta}, \qquad \tau^{\beta} = A^{\beta}_{\:\:\:\bar{\nu}}\theta^{\bar{\nu}}, \qquad A^{\alpha\beta} = A^{\beta\alpha}
\end{equation}
for a suitably determined torsion matrix $(A^{\beta}_{\:\:\:\bar{\nu}})$, where the last symmetry relation holds automatically (see \cite{W78}).  We also recall the fact that the coframe $(\theta, \theta^{1}, \ldots, \theta^{n})$ is admissible if and only if $d\theta = i g_{\alpha \bar{\beta}}\theta^{\alpha}\wedge \theta^{\bar{\beta}}$.

Now let $\hat{M}$ be another Levi-nondegenerate CR manifold of dimension $2\hat{n} + 1$, with rank $\hat{n}$ CR bundle $\hat{\mathcal{V}}$ and signature $\hat{l} \leq \hat{n}/2$.  Let $f:M \rightarrow \hat{M}$ be a smooth CR mapping in a small neighborhood of $p_{0}$.  Since our arguments are local in nature, we denote this neighborhood by $M$ also.  We use a $\hat{}$ to denote objects associated to $\hat{M}$.  Capital Latin indices $A, B$, etc. will belong to the set $\{1,\ldots,\hat{n}\}$, Greek indices $\alpha, \beta$, etc. will belong to $\{1,\ldots, n\}$, and small Latin indices $a, b$, etc. run over the complementary set $\{n+1,\ldots,\hat{n}\}$.  Let $(\theta, \theta^{\alpha})$ and $(\hat{\theta}, \hat{\theta}^{A})$ be coframes on $M$ and $\hat{M}$ respectively, and recall that $f$ is a CR mapping if
\begin{equation*}
f^{*}(\hat{\theta}) = a\theta, \qquad f^{*}(\hat{\theta}^{A}) = E^{A}_{\:\:\:\alpha}\theta^{\alpha} + E^{A}\theta,
\end{equation*}
where $a$ is a real-valued function and $E^{A}_{\:\:\:\alpha}, E^{A}$ are complex-valued functions defined near $p_{0}$.  We shall assume that $f$ is \emph{CR transversal} to $\hat{M}$ at $p_{0}$, which in our context can be expressed by saying $a(p_{0}) \neq 0$.  By applying $f^{*}$ to the equation $d\hat{\theta} = ig_{A\bar{B}}\hat{\theta}^{A}\wedge\hat{\theta}^{\bar{B}} + \hat{\theta} \wedge \phi$, we see that CR transversality of $f$ implies that $g_{\alpha\bar{\beta}} = \frac{1}{a}\hat{g}_{A\bar{B}}E^{A}_{\:\:\:\alpha}E^{\bar{B}}_{\:\:\:\beta}$.  This implies that $n \leq \hat{n}$ and $f$ is locally an embedding.  

Now suppose $(\theta, \theta^{\alpha})$ is a coframe on $M$ such that the matrix of the Levi form with respect to this coframe has $l$ negative and $n-l$ positive eigenvalues.  Let $(\hat{\theta}, \hat{\theta}^{A})$ be a coframe on $\hat{M}$ such that the matrix of the Levi form with respect to this coframe has $l'$ negative and $\hat{n} - l'$ positive eigenvalues.  If $l < n/2$ and $l' < \hat{n}/2$, we define $f$ to be \emph{side preserving} if the nonvanishing function $a$ such that $f^{*}(\hat{\theta}) = a\theta$ is strictly positive, and \emph{side reversing} if $a$ is strictly negative.  Note that this definition does not depend on the choice of pseudohermitian structure.

We state the following result, which is essentially Proposition 3.1 in \cite{BEH1} although we have been careful to distinguish the side preserving and side reversing cases.  \newtheorem{Proposition 3.1 from BEH}{Proposition}[section]
\begin{Proposition 3.1 from BEH}
Let $M$ and $\hat{M}$ be Levi-nondegenerate CR-manifolds of dimensions $2n+1$ and $2\hat{n} + 1$, and signatures $l \leq n/2$ and $l' \leq \hat{n}/2$ respectively.  Let $f:M \rightarrow \hat{M}$ be a CR mapping that is CR transversal to $\hat{M}$ along $M$.  If $(\theta, \theta^{\alpha})$ is any admissible coframe on $M$, then in a neighborhood of any point $\hat{p} \in f(M)$ in $\hat{M}$ there exists an admissible coframe $(\hat{\theta}, \hat{\theta}^{A})$ on $\hat{M}$ with $f^{*}(\hat{\theta},\hat{\theta}^{\alpha},\hat{\theta}^{a}) = (\theta, \theta^{\alpha}, 0)$.  If the Levi form of $M$ with respect to $(\theta, \theta^{\alpha})$ is constant and diagonal with $-1,\ldots,-1$ ($l$ times) and $1,\ldots,1$ ($n-l$ times) on the diagonal, then $(\hat{\theta}, \hat{\theta}^{A})$ can be chosen such that the Levi form of $\hat{M}$ relative to this coframe is constant and diagonal and if $f$ is 
\begin{description}

\item[side preserving or $l = n/2$ or $l' = \hat{n}/2$,]the diagonal elements are $-1,\ldots,-1$ ($l$ times), $1, \ldots, 1$ ($n-l$ times), $-1,\ldots, -1$ ($l' - l$ times) and $1, \ldots, 1$ ($\hat{n} - n - l' + l$ times).  With this additional property, the coframe $(\hat{\theta}, \hat{\theta}^{A})$ is uniquely determined along $M$ up to unitary transformations in $U(n, l) \times U(\hat{n} - n, l' - l)$.

\item[side reversing,]the diagonal elements are  $-1,\ldots,-1$ ($l$ times), $1, \ldots, 1$ ($n-l$ times), $-1,\ldots, -1$ ($\hat{n} - l' - l$ times) and $1, \ldots, 1$ ($l' - (n - l)$ times).  With this additional property, the coframe $(\hat{\theta}, \hat{\theta}^{A})$ is uniquely determined along $M$ up to unitary transformations in $U(n, l) \times U(\hat{n} - n, \hat{n} - l' - l)$.

\end{description}
\end{Proposition 3.1 from BEH}

Observe that if $l = n/2$ we may change the sign of $\theta$ so that the Levi form resembles the side preserving case.  If $l' = \hat{n}/2$, the two conclusions of the proposition coincide. If we fix an admissible coframe $(\theta, \theta^{\alpha})$ on $M$ and let $(\hat{\theta}, \hat{\theta}^{A})$ be an admissible coframe on $\hat{M}$ near a point $\hat{p} \in f(M)$, we shall say $(\hat{\theta}, \hat{\theta}^{A})$ is \emph{adapted} to $(\theta, \theta^{\alpha})$ on $M$ (or just to $M$ if the coframe on $M$ is understood) if it satisfies the conclusions of Proposition 3.1 with the requirement there for the Levi form.  However we will continue to write the Levi forms as $g_{\alpha\bar{\beta}}, \hat{g}_{A\bar{B}}$.  We shall also omit the $\hat{}$ over frames and coframes if there is no ambiguity.  It will be clear from the context if a form is pulled back to $M$ or not.  Under the above assumptions, we identify $M$ with the submanifold $f(M)$ and write $M \subset \hat{M}$.

Equation (3.2) implies that when $(\theta, \theta^{A})$ is adapted to $M$, if the pseudoconformal connection matrix of $(\hat{M}, \hat{\theta})$ is $\hat{\omega}_{B}^{\:\:\:A}$, then that of $(M, \theta)$ is the pullback of $\hat{\omega}_{\beta}^{\:\:\:\alpha}$.  The pulled back torsion $\hat{\tau}^{\alpha}$ is $\tau^{\alpha}$, so omitting the $\hat{}$ over these pullbacks will not cause any ambiguity and we shall do that from now on.  By the normalization of the Levi form, the second equation in (3.1) reduces to  
\begin{equation}
\omega_{B\bar{A}} + \omega_{\bar{A}B} = 0,
\end{equation}
where as before $\omega_{\bar{A}B} = \overline{\omega_{A\bar{B}}}$.

The matrix of 1-forms $(\omega_{\alpha}^{\:\:\:b})$ pulled back to $M$ defines the \emph{second fundamental form} of the embedding $f:M \rightarrow \hat{M}$.  Since $\theta^{b} = 0$ on $M$, equation (3.2) implies that on $M$, 
\begin{equation}
\omega_{\alpha}^{\:\:\:b} \wedge\theta^{\alpha} + \tau^{b}\wedge\theta = 0,
\end{equation}
and this implies that
\begin{equation}
\omega_{\alpha}^{\:\:\:b} = \omega_{\alpha \:\:\: \beta}^{\:\:\:b}\theta^{\beta}, \qquad \omega_{\alpha \:\:\: \beta}^{\:\:\:b} = \omega_{\beta \:\:\: \alpha}^{\:\:\:b}, \qquad \tau^{b} = 0.
\end{equation}

Following \cite{EHZ} we identify the CR-normal space $T_{p}^{1,0}\hat{M}/T_{p}^{1,0}M$, also denoted by $N_{p}^{1,0}\hat{M}$ with $\mathbb{C}^{\hat{n}-n}$ by choosing the equivalence classes of $L_{a}$ as a basis.  Therefore for fixed $\alpha, \beta = 1, \ldots, n$, we view the component vector $(\omega_{\alpha \:\: \beta}^{\:\: a})_{a=n+1,\ldots,\hat{n}}$ as an element of $\mathbb{C}^{\hat{n}-n}$.  By also viewing the second fundamental form as a section over $M$ of the bundle 
$T^{1,0}M\otimes N^{1,0}\hat{M} \otimes T^{1,0}M$, we may use the pseudohermitian connections on $M$ and $\hat{M}$ to define the covariant differential 
\begin{equation*}
\nabla \omega_{\alpha\:\:\beta}^{\:\:a} = d\omega_{\alpha\:\:\beta}^{\:\:a} - \omega_{\mu\:\:\beta}^{\:\:a}\omega_{\alpha}^{\:\:\mu} + \omega_{\alpha\:\:\beta}^{\:\:b}\omega_{b}^{\:\:a} - \omega_{\alpha\:\:\mu}^{\:\:a}\omega_{\beta}^{\:\:\mu}.
\end{equation*}
We write $\omega_{\alpha\:\:\beta ; \gamma}^{\:\:a}$ to denote the component in the direction $\theta^{\gamma}$ and define higher order derivatives inductively as:
\begin{equation*}
\nabla \omega_{\gamma_{1}\:\:\gamma_{2};\gamma_{3}\ldots\gamma_{j}}^{\:\:a} = d\omega_{\gamma_{1}\:\:\gamma_{2};\gamma_{3}\ldots\gamma_{j}}^{\:\:a} + \omega_{\gamma_{1}\:\:\gamma_{2};\gamma_{3}\ldots\gamma_{j}}^{\:\:b}\omega_{b}^{\:\:a} - \sum_{l=1}^{j}\omega_{\gamma_{1}\:\:\gamma_{2};\gamma_{3}\ldots\gamma_{l-1}\mu\gamma_{l+1}\ldots\gamma_{j}}^{\:\:a}\omega_{\gamma_{l}}^{\:\:\mu}.
\end{equation*}
We also consider the component vectors of higher order derivatives as elements of $\mathbb{C}^{\hat{n}-n}$ and define an increasing sequence of vector spaces
\begin{equation*}
E_{2}(p) \subseteq \ldots \subseteq E_{k}(p) \subseteq \ldots \subseteq \mathbb{C}^{\hat{n}-n}
\end{equation*}
by letting $E_{k}(p)$ be the span of the vectors
\begin{equation*}
(\omega_{\gamma_{1}\:\:\gamma_{2};\gamma_{3}\ldots\gamma_{j}}^{\:\:a})_{a=n+1,\ldots,\hat{n}}, \qquad \forall 2 \leq j \leq k, \gamma_{j}\in \{1,\ldots,n\},
\end{equation*}
evaluated at $p \in M$.  Following Lamel \cite{La01} and \cite{EHZ}, we say that the mapping $f:M \rightarrow \hat{M}$ is 
\emph{constantly $(k,s)$-degenerate} at $p$ if the vector space $E_{k}(p)$ has constant dimension $\hat{n} - n - s$ for $q$ near $p$, $E_{k+1}(q) = E_{k}(q)$, and $k$ is the smallest such integer.

\section{The Pseudoconformal Connection and Adapted $Q$-frames}

We will need the pseudoconformal connection and structure equations introduced by Chern and Moser in \cite{CM}.  Let $Y$ be the bundle of coframes $(\omega,\omega^{\alpha},\omega^{\bar{\alpha}}, \phi)$ on the real ray bundle $\pi_{E}:E\rightarrow M$ of all contact forms defining the same orientation of $M$, such that $d\omega = ig_{\alpha\bar{\beta}}\omega^{\alpha}\wedge\omega^{\bar{\beta}} + \omega\wedge\phi$ where $\omega^{\alpha} \in \pi_{E}^{*}(T'M)$ and $\omega$ is the canonical 1-form on $E$.  In \cite{CM} it was shown that these forms can be completed to a full set of invariants on $Y$ given by the coframe of 1-forms
\[ (\omega, \omega^{\alpha}, \omega^{\bar{\alpha}}, \phi, \phi_{\beta}^{\alpha}, \phi^{\alpha}, \phi^{\bar{\alpha}}, \psi) \] which define the pseudoconformal connection on $Y$.  These forms satisfy the structure equations, which we will use extensively (see \cite{CM} and its appendix):
\begin{align}
& \phi_{\alpha\bar{\beta}} + \phi_{\bar{\beta}\alpha} = g_{\alpha\bar{\beta}}\phi, \nonumber \\
& d\omega = i\omega^{\mu}\wedge\omega_{\mu} + \omega \wedge \phi, \nonumber \\
& d\omega^{\alpha} = \omega^{\mu}\wedge\phi_{\mu}^{\:\:\alpha} + \omega\wedge\phi^{\alpha}, \nonumber \\
& d\phi = i\omega_{\bar{\nu}}\wedge\phi^{\bar{\nu}} + i\phi_{\bar{\nu}}\wedge\omega^{\bar{\nu}} + \omega \wedge \psi, \nonumber \\
& d\phi_{\beta}^{\:\:\alpha} = \phi_{\beta}^{\:\:\mu}\wedge\phi_{\mu}^{\:\:\alpha} + i\omega_{\beta}\wedge\phi^{\alpha} - i\phi_{\beta}\wedge\omega^{\alpha} -i\delta_{\beta}^{\:\:\alpha}\phi_{\mu}\wedge\omega^{\mu} - \frac{\delta_{\beta}^{\:\:\alpha}}{2}\psi\wedge\omega + \Phi_{\beta}^{\:\:\alpha}, \nonumber \\
& d\phi^{\alpha} = \phi\wedge\phi^{\alpha} + \phi^{\mu}\wedge\phi_{\mu}^{\:\:\alpha} - \frac{1}{2}\psi\wedge\omega^{\alpha} + \Phi^{\alpha}, \nonumber \\
& d\psi = \phi\wedge\psi + 2i\phi^{\mu}\wedge\phi_{\mu} + \Psi. 
\end{align}

Here the 2-forms $\Phi_{\beta}^{\:\:\alpha}, \Phi^{\alpha}, \Psi$ give the \emph{pseudoconformal curvature} of $M$.  We may decompose $\Phi_{\beta}^{\:\:\alpha}$ as follows
\begin{equation*}
\Phi_{\beta}^{\:\:\alpha} = S_{\beta\:\:\:\mu\bar{\nu}}^{\:\:\alpha}\omega^{\mu}\wedge\omega^{\bar{\nu}} + V_{\beta\:\:\:\mu}^{\:\:\alpha}\omega^{\mu}\wedge\omega + V^{\alpha}_{\:\:\:\beta\bar{\nu}}\omega\wedge\omega^{\bar{\nu}}.
\end{equation*}
We will also refer to the tensor $S_{\beta\:\:\:\mu\bar{\nu}}^{\:\:\alpha}$ as the pseudoconformal curvature of $M$.  We require $ S_{\beta\:\:\:\mu\bar{\nu}}^{\:\:\alpha}$ to satisfy certain trace and symmetry conditions (see \cite{CM}), but for the purposes of this paper, the important point to emphasize is that \emph{for a hyperquadric, the pseudoconformal curvature vanishes}.

If we fix a contact form $\theta$ that defines a section $M \rightarrow E$, then any admissible coframe $(\theta, \theta^{\alpha})$ for $M$ defines a unique section $M \rightarrow Y$ under which the pullbacks of $(\omega, \omega^{\alpha})$ coincide with $(\theta, \theta^{\alpha})$ and the pullback of $\phi$ vanishes.  As in \cite{W78} we use this section to pull the pseudoconformal connection forms back to $M$.  Although the pulled back forms on $M$ now depend on the choice of admissible coframe, we shall use the same notation, and thus we have
\begin{equation*}
\theta = \omega,\qquad \theta^{\alpha} = \omega^{\alpha},\qquad \phi = 0
\end{equation*}
on $M$.
As in \cite{W78}, we may write the pulled back tangential pseudoconformal curvature tensor $ S_{\beta\:\:\:\mu\bar{\nu}}^{\:\:\alpha}$ in terms of the tangential pseudohermitian curvature tensor $ R_{\beta\:\:\:\mu\bar{\nu}}^{\:\:\alpha}$ by
\begin{equation*}
 S_{\alpha\bar{\beta}\mu\bar{\nu}} = R_{\alpha\bar{\beta}\mu\bar{\nu}} - \frac{R_{\alpha\bar{\beta}}g_{\mu\bar{\nu}} + R_{\mu\bar{\beta}}g_{\alpha\bar{\nu}} + R_{\alpha\bar{\nu}}g_{\mu\bar{\beta}} + R_{\mu\bar{\nu}}g_{\alpha\bar{\beta}}}{n+2} + \frac{R(g_{\alpha\bar{\beta}}g_{\mu\bar{\nu}} + g_{\alpha\bar{\nu}}g_{\mu\bar{\beta}})}{(n+1)(n+2)},
\end{equation*}
where
\begin{equation*}
R_{\alpha\bar{\beta}} := R_{\mu\:\:\:\alpha\bar{\beta}}^{\:\:\mu} \qquad \textrm{and } R:= R_{\mu}^{\:\:\mu}
\end{equation*}
are respectively the pseudohermitian Ricci and scalar curvature of $(M, \theta)$.  This formula expresses the fact that $S_{\alpha\bar{\beta}\mu\bar{\nu}}$ is the ``traceless component" of $R_{\alpha\bar{\beta}\mu\bar{\nu}}$ with respect to the decomposition of the space of all tensors with the symmetry conditions of $S_{\alpha\bar{\beta}\mu\bar{\nu}}$ into the direct sum of the subspace of tensors with trace zero and the subspace of \emph{conformally flat tensors}, i.e. tensors of the form
\begin{equation}
T_{\alpha\bar{\beta}\mu\bar{\nu}} = H_{\alpha\bar{\beta}}g_{\mu\bar{\nu}} + H_{\mu\bar{\beta}}g_{\alpha\bar{\nu}} + H_{\alpha\bar{\nu}}g_{\mu\bar{\beta}} + H_{\mu\bar{\nu}}g_{\alpha\bar{\beta}},
\end{equation}
where $(H_{\alpha\bar{\beta}})$ is any Hermitian matrix.  We shall call two tensors as above \emph{conformally equivalent} if their difference is of the form of equation (4.2).  Note that covariant derivatives of conformally flat tensors are conformally flat, because $\nabla g_{\alpha\bar{\beta}} = 0$.

The following result relates the pseudoconformal and pseudohermitian connection forms.  It is alluded to in \cite{W78} and a proof may be found in \cite{EHZ}, where the result appears as Proposition 3.1.  Note that although the Proposition in \cite{EHZ} is stated only for $M$ strictly pseudoconvex, the result is valid in the Levi-nondegenerate situation.
\newtheorem{Proposition 3.2}[Proposition 3.1 from BEH]{Proposition}
\begin{Proposition 3.2}
Let $M$ be a smooth Levi-nondegenerate CR-manifold of hypersurface type with CR dimension $n$, and with respect to an admissible coframe $(\theta, \theta^{\alpha})$ let the pseudoconformal and pseudohermitian connection forms be pulled back to $M$ as above.  Then we have the following relations:
\begin{equation*}
\phi_{\beta}^{\:\:\alpha} = \omega_{\beta}^{\:\:\alpha} + D_{\beta}^{\:\:\alpha}\theta, \qquad \phi^{\alpha} = \tau^{\alpha} + D_{\mu}^{\:\:\alpha}\theta^{\mu} + E^{\alpha}\theta, \qquad \psi = iE_{\mu}\theta^{\mu} - iE_{\bar{\nu}}\theta^{\bar{\nu}} + B\theta,
\end{equation*}
where
\begin{align*}
&D_{\alpha\bar{\beta}} := \frac{iR_{\alpha\bar{\beta}}}{n+2} - \frac{iRg_{\alpha\bar{\beta}}}{2(n+1)(n+2)}, \nonumber \\
&E^{\alpha} := \frac{2i}{2n + 1}(A^{\alpha\mu}_{\:\:\:\:;\mu} - D^{\bar{\nu}\alpha}_{\:\:\:\:;\bar{\nu}}), \nonumber \\
&B := \frac{1}{n}(E^{\mu}_{\:\:\:;\mu} + E^{\bar{\nu}}_{\:\:\:;\bar{\nu}} - 2A^{\beta\mu}A_{\beta\mu} + 2D^{\bar{\nu}\alpha}D_{\bar{\nu}\alpha}).
\end{align*}
\end{Proposition 3.2}

Another notion that will prove useful is that of an adapted $Q$-frame.  We embed $\mathbb{C}^{\hat{n}+1}$ in $\mathbb{CP}^{\hat{n}+1}$ as the set $\{\zeta^{0} \neq 0 \}$ in the homogeneous coordinates $[\zeta^{0}:\zeta^{1}:\ldots:\zeta^{\hat{n}+1}]$, and following section 1 of \cite{CM}, realize the the quadric $Q_{l}^{\hat{n}}$ given in $\mathbb{CP}^{\hat{n}+1}$ by the equation $(\zeta, \zeta) = 0$, where the Hermitian scalar product $(\cdot, \cdot)$ is defined by
\begin{equation}
(\zeta, \tau) := \hat{I}_{A\bar{B}}\zeta^{A}\overline{\tau^{B}} + \frac{i}{2}(\zeta^{\hat{n} + 1}\overline{\tau^{0}} - i \zeta^{0}\overline{\tau^{\hat{n}+1}}).
\end{equation}
In the above, $\hat{I}$ is the diagonal matrix with first $l$ diagonal entries equal to $-1$ and all subsequent diagonal entries equal $1$.  A \emph{$Q$-frame} (see e.g. \cite{CM}) is a unimodular basis $(Z_{0}, \ldots, Z_{\hat{n}+1})$ of $\mathbb{C}^{\hat{n}+2}$ such that $Z_{0}$ and $Z_{\hat{n}+1}$, as points in $\mathbb{CP}^{\hat{n}+1}$, are on $Q$, the vectors $(Z_{A})$ form an orthonormal basis (relative to the inner product (4.3)) for the complex tangent space to the quadric at $Z_{0}$ and $Z_{\hat{n}+1}$, and $(Z_{\hat{n}+1}, Z_{0}) = i/2$.  We will denote the corresponding points in $\mathbb{CP}^{\hat{n}+1}$ also by $Z_{0}$ and $Z_{\hat{n}+1}$; it should be clear from the context whether the point is in $\mathbb{C}^{\hat{n}+2}$ or $\mathbb{CP}^{\hat{n}+1}$ .

On the space $\mathfrak{B}$ of all $Q$-frames there is a natural free transitive action of the group $\textbf{SU}(l+1, \hat{n} - l +1)$ of unimodular $(\hat{n}+2)\times(\hat{n}+2)$ matrices that preserve the inner product (4.3).  Hence, any fixed $Q$-frame defines an isomorphism between $\mathfrak{B}$ and $\textbf{SU}(l+1, \hat{n} - l +1)$.  On the space $\mathfrak{B}$, there are \emph{Maurer-Cartan forms} $\pi_{\Lambda}^{\:\:\:\Omega}$, where capital Greek indices run from $0$ to $\hat{n}+1$, defined by 
\begin{equation}
dZ_{\Lambda} = \pi_{\Lambda}^{\:\:\:\Omega}Z_{\Omega}
\end{equation}
and satisfying $d\pi_{\Lambda}^{\:\:\:\Omega} = \pi_{\Lambda}^{\:\:\:\Gamma} \wedge \pi_{\Gamma}^{\:\:\:\Omega}$.  Here the natural $\mathbb{C}^{\hat{n}+2}$ valued 1-forms $dZ_{\Lambda}$ on $\mathfrak{B}$ are defined as differentials of the map $(Z_{0}, \ldots, Z_{\hat{n}+1}) \rightarrow Z_{\Lambda}$.

Recall from \cite{CM} and \cite{W79} that a smoothly varying $Q$-frame $(Z_{\Lambda}) = (Z_{\Lambda}(p))$ for $p \in Q$ is said to be \emph{adapted} to $Q$ if $Z_{0}(p) = p$ as points of $\mathbb{CP}^{\hat{n}+1}$.  It is shown in section 5 of \cite{CM} that if we use an adapted $Q$-frame to pull back the 1-forms $\pi_{\Lambda}^{\:\:\:\Omega}$ from $\mathfrak{B}$ to $Q$ and set
\begin{equation}
\theta := \frac{1}{2}\pi_{0}^{\:\:\:\hat{n}+1}, \qquad \theta^{A} := \pi_{0}^{\:\:\:A}, \qquad \xi := -\pi_{0}^{\:\:\:0} + \overline{\pi_{0}^{\:\:\:0}},
\end{equation}
we obtain a coframe $(\theta, \theta^{A})$ on $Q$ and a form $\xi$ satisfying the structure equation
\begin{equation*}
d\theta = i \hat{I}_{A\bar{B}}\theta^{A}\wedge\theta^{\bar{B}} + \theta \wedge \xi.
\end{equation*}
In particular, it follows from (4.5) that the coframe $(\theta^{A}, 2\theta)$ is dual to the frame defined by $(Z_{A}, Z_{\hat{n}+1})$ on $Q$ and hence depends only on the values of $(Z_{\Lambda})$ at the same points.  Furthermore, there is a unique section $M \rightarrow Y$ for which the pullbacks of the forms $(\omega, \omega^{\alpha}, \phi)$ are $(\theta, \theta^{\alpha}, \xi)$ respectively.  Then the pulled back forms $(\hat{\phi}_{B}^{\:\:\:A}, \hat{\phi}^{A}, \hat{\psi})$ are given by  (5.8b) from \cite{CM}:
\begin{equation}
\hat{\phi}_{B}^{\:\:\:A} = \pi_{B}^{\:\:\:A} - \delta_{B}^{\:\:\:A}\pi_{0}^{\:\:\:0}, \qquad \hat{\phi}^{A} = 2\pi_{\hat{n}+1}^{\:\:\:A}, \qquad \hat{\psi} = -4\pi_{\hat{n}+1}^{\:\:\:0}.
\end{equation}
As in (5.30) from \cite{CM}, the pulled back forms $\pi_{\Lambda}^{\:\:\:\Omega}$ can be uniquely solved from (4.5-4.6):
\begin{align}
(\hat{n} + 2)\pi_{0}^{\:\:\:0} &= -\hat{\phi}_{C}^{\:\:\:C} - \xi, & \pi_{0}^{\:\:\:A} &= \theta^{A}, \qquad & \pi_{0}^{\:\:\:\hat{n}+1} &= 2\theta, \nonumber \\
\pi_{A}^{\:\:\:0} &= -i\hat{\phi}_{A}, & \pi_{A}^{\:\:\:B}  &= \hat{\phi}_{A}^{\:\:\:B} + \delta_{A}^{\:\:\:B}\pi_{0}^{\:\:\:0}, \qquad & \pi_{A}^{\:\:\:\hat{n}+1} &= 2i\theta_{A}, \nonumber \\
4\pi_{\hat{n}+1}^{\:\:\:0} &= -\hat{\psi}, \qquad & 2\pi_{\hat{n}+1}^{\:\:\:A} &= \hat{\phi}^{A}, \qquad & (\hat{n}+2)\pi_{\hat{n}+1}^{\:\:\:\hat{n}+1} &= \hat{\phi}_{\bar{D}}^{\:\:\:\bar{D}} + \xi.
\end{align}
Thus, the pullback of $\pi_{\Lambda}^{\:\:\:\Omega}$ is completely determined by the pullbacks $(\theta, \theta^{A}, \xi, \hat{\phi}_{B}^{\:\:\:A}, \hat{\phi}^{A}, \hat{\psi})$.  Following section 8 of \cite{EHZ}, we note that for any choice of an admissible coframe $(\theta, \theta^{A})$ on $Q$, there exists an adapted $Q$-frame $(Z_{\Lambda})$ such that (4.7) holds with $\xi = 0$.

\section{Proof of Theorem 1.1}

The following lemma will be a key ingredient in the proof of Theorem 1.1:
\newtheorem{Key Lemma}{Lemma}[section]
\begin{Key Lemma}

Let $g$ be a diagonal matrix in $\mathbb{C}^{d}$ with either positive or negative $1$ in each diagonal entry and denote by $e_{j} = (0, \ldots, 1, \ldots, 0)^{T}$ the $j^{th}$ standard basis vector in $\mathbb{C}^{d}$.  Let $E$ be the span of $r$ independent vectors in $\mathbb{C}^{d}$, with $r+s = d$.  Without loss of generality, suppose $E$ is a graph over $\{e_{s+1}, \ldots, e_{d}\}$, that is, there exists a $d\times r$ matrix of the form $\begin{pmatrix} C^{T} \\ I \end{pmatrix}$ where $C^{T}$ is $s\times r$, whose columns span $E$.  Then there exists an invertible matrix $A$ in $\mathbb{C}^{d}$ such that if $N = A^{-1}$, then for $v \in E$, $N^{T}v \in \textrm{span}\{e_{s+1}, \ldots, e_{d}\}$ and if $\tilde{g}:= A^{*}gA$, then $\tilde{g}_{pq} = 0$ when $p \in \{s+1,\ldots,d\}$ and $q \in \{1, \ldots, s\}$.

\begin{proof}
Define $I_{1}$ and $I_{2}$ to be the $s\times s$ and $r\times r$ upper left and lower right blocks of $g$, respectively.  Choose a matrix norm such that $||I_{j}|| \le 1$ for $j=1,2$ and nonzero constant $\lambda$ such that $|\lambda|^{2} > \max\{||C^{*}I_{2}C||, ||I_{2}CI_{1}C^{*}I_{2}||\}$.

We now show that $A := \begin{pmatrix}\lambda I && -\frac{1}{\bar{\lambda}}I_{1}C^{*}I_{2} \\ C && I \end{pmatrix}$, where the upper left block is $s\times s$ and the lower right block is $r \times r$ satisfies the desired requirements.  Note that by construction, $A^{T}$ carries the span of $\{e_{s+1}, \ldots, e_{d}\}$ to $E$, so $N^{T}$ takes $E$ to the span of $\{e_{s+1}, \ldots, e_{d}\}$.

We compute $A^{*}gA:$
\begin{align*} A^{*}gA &= \begin{pmatrix}\bar{\lambda} I && C^{*} \\ -\frac{1}{\lambda}I_{2}CI_{1} && I \end{pmatrix}\begin{pmatrix}I_1 && 0\\ 0 && I_2 \end{pmatrix}\begin{pmatrix}\lambda I && -\frac{1}{\bar{\lambda}}I_{1}C^{*}I_{2} \\ C && I \end{pmatrix} \nonumber \\
&= \begin{pmatrix}\bar{\lambda} I && C^{*} \\ -\frac{1}{\lambda}I_{2}CI_{1} && I \end{pmatrix}\begin{pmatrix}\lambda I_{1} && -\frac{1}{\bar{\lambda}}C^{*}I_{2} \\ I_{2}C && I_{2} \end{pmatrix} \nonumber \\
&= \begin{pmatrix}|\lambda|^{2}(I_{1} + \frac{1}{|\lambda|^{2}}C^{*}I_{2}C) && 0 \\ 0 && I_{2} + \frac{1}{|\lambda|^{2}}I_{2}CI_{1}C^{*}I_{2} \end{pmatrix}.
\end{align*}
This shows that $A^{*}gA$ is block diagonal.  To see that $A$ is invertible, it suffices to show each block of $A^{*}gA$ is invertible.  Up to a constant, each block is of the form $I_{j} + L$, where $L$ has norm less than $1$ by our choice of $\lambda$.  This implies that $I + I_{j}L$ is invertible (with the appropriate dimensions of $I$ in each block), so there is a matrix $D$ such that $(I + I_{j}L)D = I$.  Hence by multiplying both sides on the left and right by $I_{j}$ we have $(I_{j} + L)DI_{j} = I$, so $I_{j} + L$ is invertible, as desired.
\end{proof}
\end{Key Lemma}

\emph{Proof of Theorem 1.1.}  We choose an admissible coframe $(\theta, \theta^{A})$ on $Q$ near $f(p)$ adapted to an admissible coframe $(\theta, \theta^{\alpha})$ on $M$ and denote by $(\omega_{\alpha\textrm{ }\beta}^{\phantom{1}a})$ the second fundamental form of $f$ relative to this coframe.  Since the mapping $f$ is $(k,s)$-degenerate near $p$, we have that the dimension of span$\{\omega_{\gamma_{1} \phantom{1} \gamma_{2};\gamma_{3}\ldots\gamma_{t}}^{\phantom{1}a}, 2\le t \le k\}$ is $r = d-s$ near $p$.  We introduce some notation; the indices $\ast, \#$ run over the set $n+1, \ldots, n+r$ (possibly empty) and the indices $i, j$ run over the set $n + r + 1, \ldots, N$.

We now  fix $\alpha, \beta$ and identify $(\omega_{\alpha\textrm{ }\beta}^{\phantom{1}a}(p))$ as a vector in $\mathbb{C}^{N-n}$.  We apply Lemma 5.1 with $g_{a\bar{b}}$ as the matrix $g$ and after the above identification, we let $E =$ span$\{\omega_{\gamma_{1} \phantom{1} \gamma_{2};\gamma_{3}\ldots\gamma_{t}}^{\phantom{1}a}, 2\le t \le k\}$.  This produces a smooth matrix-valued function $A$.  We change basis (only on the normal space) via
$\begin{pmatrix}\theta^{n+1} \\ \vdots \\ \theta^{N} \end{pmatrix} = \begin{pmatrix} && && \\ && A && \\ && && \end{pmatrix} \begin{pmatrix} \tilde{\theta}^{n+1} \\ \vdots \\ \tilde{\theta}^{N} \end{pmatrix}$
then we have
\begin{equation} \textrm{span}\{\omega_{\gamma_{1} \phantom{1} \gamma_{2};\gamma_{3}\ldots\gamma_{t}}^{\phantom{1}\#}L_{\#}, 2\le t \le k\} = \textrm{span}\{\tilde{L}_{\#}\}, \textrm{   and} \:\:\: \omega_{\gamma_{1} \phantom{1} \gamma_{2};\gamma_{3}\ldots\gamma_{t}}^{\phantom{1}j} \equiv 0, \phantom{1} t \ge 2. \end{equation}
We now relabel and omit the tilde notation.  Note that our Levi form on the normal space is no longer necessarily constant, but does satisfy at each point the conclusion of Lemma 5.1, so $g_{\# \bar{j}} = 0$.  Also, we still have the relations $f^{*}(\theta^{a})=0$ and $\hat{g}_{\alpha \bar{\beta}} = g_{\alpha \bar{\beta}}.$  Note that the inverse of a block diagonal matrix is block diagonal, so $g^{A\bar{B}}$ has the same form as $g_{A\bar{B}}.$

Because $\hat{\omega}_{\#}^{\phantom{2}j}$ is a 1-form on $M$, we have
\begin{equation} \hat{\omega}_{\#}^{\phantom{2}j} = \hat{\omega}_{\# \phantom{2}\mu}^{\phantom{2}j}\theta^{\mu} + \hat{\omega}_{\# \phantom{2}\bar{\nu}}^{\phantom{2}j}\theta^{\bar{\nu}} + \hat{\omega}_{\# \phantom{2} 0 }^{\phantom{2}j}\theta
\end{equation}
for suitable coefficients.

Now by the definition of covariant derivative, we have
\begin{align}
\nabla \omega_{\gamma_{1} \phantom{1} \gamma_{2};\gamma_{3}\ldots\gamma_{t}}^{\phantom{2}j} &= d\omega_{\gamma_{1} \phantom{1} \gamma_{2};\gamma_{3}\ldots\gamma_{t}}^{\phantom{2}j} + \omega_{\gamma_{1} \phantom{1} \gamma_{2};\gamma_{3}\ldots\gamma_{t}}^{\phantom{2}i}\hat{\omega}_{i}^{\phantom{2}j} + \omega_{\gamma_{1} \phantom{1} \gamma_{2};\gamma_{3}\ldots\gamma_{t}}^{\phantom{2}\#}\hat{\omega}_{\#}^{\phantom{2}j}  \nonumber \\
&\qquad - \sum_{q=1}^{t}\omega_{\gamma_{1} \phantom{1} \gamma_{2};\gamma_{3}\ldots\gamma_{q-1}\mu\gamma_{q+1}\ldots\gamma_{t}}^{\phantom{2}j}\hat{\omega}_{\gamma_{q}}^{\phantom{2}\mu} \nonumber
\end{align}
so by $(5.1)$ we have
\[ \nabla \omega_{\gamma_{1} \phantom{1} \gamma_{2};\gamma_{3}\ldots\gamma_{t}}^{\phantom{2}j} = \omega_{\gamma_{1} \phantom{1} \gamma_{2};\gamma_{3}\ldots\gamma_{t}}^{\phantom{2}\#}\hat{\omega}_{\#}^{\phantom{2}j}. \]
This implies that
\[\omega_{\gamma_{1} \phantom{1} \gamma_{2};\gamma_{3}\ldots\gamma_{t}\mu}^{\phantom{2}j} = \omega_{\gamma_{1} \phantom{1} \gamma_{2};\gamma_{3}\ldots\gamma_{t}}^{\phantom{2}\#}\hat{\omega}_{\# \phantom{1} \mu}^{\phantom{2}j}\]
and because the left side is zero we have
\begin{equation}
\omega_{\gamma_{1} \phantom{1} \gamma_{2};\gamma_{3}\ldots\gamma_{t}}^{\phantom{2}\#}\hat{\omega}_{\# \phantom{1} \mu}^{\phantom{2}j} = 0.
\end{equation}
Now if $j, \mu$ are fixed and $\hat{\omega}_{\# \phantom{1} \mu}^{\phantom{2}j} \ne 0$ for some $\#$ then pick $r$ independent vectors with $r$ components $(\omega_{\gamma_{1} \phantom{1} \gamma_{2};\gamma_{3}\ldots\gamma_{t}}^{\phantom{2}\ast})$, make a matrix $B$ with these as the rows and let $v$ be the vector $(\hat{\omega}_{\# \phantom{1} \mu}^{\phantom{2}j})$ as $\#$ varies.  Then $Bv=0$ contradicting independence of the rows of $B$.  This implies that
\begin{equation}
\hat{\omega}_{\# \phantom{1} \mu}^{\phantom{2}j} = 0.
\end{equation}

Now applying Proposition 4.1, and noting that, by equations (3.5) and (5.1) we have $\hat{\omega}_{\alpha}^{\:\:\:j} = 0$ and $\tau^{a} = 0$, we find
\begin{equation}
\hat{\phi}_{\alpha}^{\phantom{2}j} = \hat{D}_{\alpha}^{\phantom{2}j}\theta, \qquad \hat{\phi}^{j} = \hat{D}_{\mu}^{\phantom{2}j}\theta^{\mu} + \hat{E}^{j}\theta,
\end{equation}
and
\begin{equation}
\hat{\phi}_{\alpha}^{\phantom{2}\#} = \hat{\omega}_{\alpha}^{\phantom{2}\#} + \hat{D}_{\alpha}^{\phantom{2}\#}\theta, \qquad \hat{\phi}^{\#} = \hat{D}_{\mu}^{\phantom{2}\#}\theta^{\mu} + \hat{E}^{\#}\theta.
\end{equation}

Next, we differentiate $\hat{\phi}_{\alpha}^{\phantom{2}j}$ and compute mod $\theta$ to obtain
\[d\hat{\phi}_{\alpha}^{\phantom{2}j} \equiv \hat{D}_{\alpha}^{\phantom{2}j}g_{\mu\bar{\nu}}\theta^{\mu}\wedge\theta^{\bar{\nu}} \qquad \textrm{mod }\theta\]
On the other hand, we may compute $d\hat{\phi}_{\alpha}^{\phantom{2}j}$ mod $\theta$ using the structure equations (4.1).  We have
\begin{align*}d\hat{\phi}_{\alpha}^{\phantom{2}j} &\equiv \hat{\phi}_{\alpha}^{\phantom{2}A}\wedge\hat{\phi}_{A}^{\phantom{2}j} + i\theta_{\alpha}\wedge\hat{\phi}^{j} - i\phi_{\alpha}\wedge\theta^{j} - i\delta_{\alpha}^{\phantom{2}j}\phi_{A}\wedge\theta^{A} - \frac{\delta_{\alpha}^{\phantom{2}j}}{2}\psi\wedge\theta + \Phi_{\alpha}^{\phantom{2}j} \nonumber \\
&\equiv \hat{\phi}_{\alpha}^{\phantom{2}A}\wedge\hat{\phi}_{A}^{\phantom{2}j} + i\theta_{\alpha}\wedge\hat{\phi}^{j} \qquad \textrm{mod }\theta. \nonumber
\end{align*}
We note that in the structure equation above the third term is zero because the pullback of $\theta^{j}$ vanishes, the fourth and fifth terms are zero because of the indices of the kronecker delta, and the last term is zero because of the vanishing pseudoconformal curvature of the target hyperquadric.

We expand the above to obtain
\begin{align*}d\hat{\phi}_{\alpha}^{\phantom{2}j} &\equiv \hat{\phi}_{\alpha}^{\phantom{2}\beta}\wedge\hat{\phi}_{\beta}^{\phantom{2}j} + \hat{\phi}_{\alpha}^{\phantom{2}\#}\wedge\hat{\phi}_{\#}^{\phantom{2}j} + \hat{\phi}_{\alpha}^{\phantom{2}i}\wedge\hat{\phi}_{i}^{\phantom{2}j} + i\theta_{\alpha}\wedge\hat{\phi}^{j} \nonumber \\
&\equiv \hat{\phi}_{\alpha}^{\phantom{2}\#}\wedge\hat{\phi}_{\#}^{\phantom{2}j} + i\theta_{\alpha}\wedge\hat{\phi}^{j} \nonumber \\
&\equiv \hat{\omega}_{\alpha\phantom{2}\mu}^{\phantom{2}\#}\theta^{\mu}\wedge\hat{\phi}_{\#}^{\phantom{2}j} - i\hat{\phi}^{j}\wedge g_{\alpha\bar{A}}\theta^{\bar{A}} \qquad \textrm{mod }\theta. \nonumber
\end{align*}
In the second equivalence we used equation (5.5) and computed mod $\theta$, and in the last equivalence we used both (5.5) and equation (3.5).

Now we may put these equations together and group terms to obtain
\begin{equation}
\hat{\omega}_{\alpha\phantom{2}\mu}^{\phantom{2}\#}\theta^{\mu}\wedge\hat{\phi}_{\#}^{\phantom{2}j} \equiv i(g_{\alpha\bar{\nu}}\hat{D}_{\mu}^{\phantom{2}j} + g_{\mu\bar{\nu}}\hat{D}_{\alpha}^{\phantom{2}j})\theta^{\mu}\wedge\theta^{\bar{\nu}} \qquad \textrm{mod }\theta.
\end{equation}
By Proposition 4.1 and equation (5.2), we compute $\hat{\phi}_{\#}^{\phantom{2}j}$ and identify the coefficients of $\theta^{\mu}\wedge\theta^{\bar{\nu}}$ to obtain
\[\hat{\omega}_{\alpha\phantom{2}\mu}^{\phantom{2}\#}\hat{\omega}_{\#\phantom{2}\bar{\nu}}^{\phantom{2}j} = i(g_{\alpha\bar{\nu}}\hat{D}_{\mu}^{\phantom{2}j} + g_{\mu\bar{\nu}}\hat{D}_{\alpha}^{\phantom{2}j}).\]
This holds in a neighborhood of $p$, so we now work at a point $q$ close to $p$.  Let
\[f_{\#}(z) = \hat{\omega}_{\alpha\phantom{2}\mu}^{\phantom{2}\#}z^{\alpha}z^{\mu} \qquad \textrm{ and } \qquad g_{\#}(z) = \hat{\omega}_{\bar{\#}\phantom{2}\bar{\nu}}^{\phantom{2}\bar{j}}z^{\nu}\]
where $\hat{\omega}_{\bar{\#}}^{\phantom{2}\bar{j}} = \hat{\omega}_{\bar{\#}\phantom{2}\mu}^{\phantom{2}\bar{j}}\theta^{\mu} + \hat{\omega}_{\bar{\#}\phantom{2}\bar{\nu}}^{\phantom{2}\bar{j}}\theta^{\bar{\nu}} + \hat{\omega}_{\bar{\#}\phantom{2}0}^{\phantom{2}\bar{j}}\theta.$
Then we have that
\begin{align*}
\sum_{\#}f_{\#}(z)\overline{g_{\#}(z)} &= \hat{\omega}_{\alpha\phantom{2}\mu}^{\phantom{2}\#}\hat{\omega}_{\#\phantom{2}\bar{\nu}}^{\phantom{2}j}z^{\alpha}z^{\mu}z^{\bar{\nu}} \nonumber \\
&= i(g_{\alpha\bar{\nu}}\hat{D}_{\mu}^{\phantom{2}j} + g_{\mu\bar{\nu}}\hat{D}_{\alpha}^{\phantom{2}j})z^{\alpha}z^{\mu}z^{\bar{\nu}} \nonumber \\
&= \langle z, z\rangle_{g}(i\hat{D}_{\mu}^{\phantom{2}j}z^{\mu} + i\hat{D}_{\alpha}^{\phantom{2}j}z^{\alpha}) \nonumber.
\end{align*}
Therefore by Lemma 2.1, since $\#$ runs over an index set of size $r$ and by assumption $r = N - n - s < n$, we have
\begin{equation}
\hat{\omega}_{\alpha\phantom{2}\mu}^{\phantom{2}\#}\hat{\omega}_{\#\phantom{2}\bar{\nu}}^{\phantom{2}j} = 0.
\end{equation}
This implies that $g_{\alpha\bar{\nu}}\hat{D}_{\mu}^{\phantom{2}j} + g_{\mu\bar{\nu}}\hat{D}_{\alpha}^{\phantom{2}j} = 0.$  Let $\alpha = \mu$ and choose $\bar{\nu}$ such that $g_{\alpha\bar{\nu}} \ne 0$, which exists since no row is completely zero.  This implies $\hat{D}_{\alpha}^{\phantom{2}j} = 0$, so
\begin{equation}
\hat{\phi}_{\alpha}^{\phantom{2}j} = 0, \qquad \hat{\phi}^{j} = \hat{E}^{j}\theta.
\end{equation}

Combining the structure equation for $d\hat{\phi}_{\alpha}^{\phantom{2}j}$ with the above result yields
\[0 = \hat{\phi}_{\alpha}^{\phantom{2}A}\wedge\hat{\phi}_{A}^{\phantom{2}j} + i\theta_{\alpha}\wedge\hat{\phi}^{j} - i\hat{\phi}_{\alpha}\wedge\hat{\theta}^{j}.\]
We only consider those terms containing a $\theta^{\mu}\wedge\theta$ and discover, using Proposition 4.1 and equation (5.2), that
\begin{align*} 0 &= \hat{\phi}_{\alpha}^{\phantom{2}\#}\wedge\hat{\phi}_{\#}^{\phantom{2}j} \nonumber \\
&= (\omega_{\alpha\phantom{2}\mu}^{\phantom{2}\#}\theta^{\mu} + \hat{D}_{\alpha}^{\phantom{2}\#}\theta)\wedge(\hat{\omega}_{\#\phantom{2}\mu}^{\phantom{2}j}\theta^{\mu} + \hat{\omega}_{\#\phantom{2}\bar{\nu}}^{\phantom{2}j}\theta^{\bar{\nu}} + (\hat{\omega}_{\#\phantom{2}0}^{\phantom{2}j} + \hat{D}_{\#}^{\phantom{2}j})\theta). \nonumber
\end{align*}
Keeping the $\theta^{\mu}\wedge\theta$ terms and using equation (5.4), we obtain
\begin{equation}
0 = \omega_{\alpha\phantom{2}\mu}^{\phantom{2}\#}(\hat{\omega}_{\#\phantom{2}0}^{\phantom{2}j} + \hat{D}_{\#}^{\phantom{2}j}).
\end{equation}

Now we would like to show that $\hat{\phi}_{\#}^{\phantom{2}j}=0$, so since $\hat{\phi}_{\#}^{\phantom{2}j} = \hat{\omega}_{\#\phantom{2}\bar{\nu}}^{\phantom{2}j}\theta^{\bar{\nu}} + (\hat{\omega}_{\#\phantom{2}0}^{\phantom{2}j} + \hat{D}_{\#}^{\phantom{2}j})\theta$ by Proposition 4.1 and equation (5.4), it suffices to show
\begin{equation}
\omega_{\gamma_{1}\phantom{1}\gamma_{2};\gamma_{3}\ldots\gamma_{t}}^{\phantom{9}\#}\hat{\omega}_{\#\phantom{2}\bar{\nu}}^{\phantom{2}j} = \omega_{\gamma_{1}\phantom{1}\gamma_{2};\gamma_{3}\ldots\gamma_{t}}^{\phantom{9}\#}(\hat{\omega}_{\#\phantom{2}0}^{\phantom{2}j} + \hat{D}_{\#}^{\phantom{2}j}) = 0, \qquad t\ge2
\end{equation}
by the same reason equation (5.3) implied (5.4).

Before proving (5.11), we first wish to show that $\hat{\omega}_{\#\phantom{2}\bar{\nu};\mu}^{\phantom{2}j}$ is a sum of multiples of the Levi form.  We differentiate the expression for $\hat{\phi}_{\#}^{\phantom{2}j}$ in Proposition 4.1, set it equal to the corresponding structure equation, and compute mod $\theta$ to obtain
\[\hat{\phi}_{\#}^{\phantom{2}A}\wedge\hat{\phi}_{A}^{\phantom{2}j} = d\omega_{\#}^{\phantom{2}j} + D_{\#}^{\phantom{2}j}g_{\mu\bar{\nu}}\theta^{\mu}\wedge\theta^{\bar{\nu}}.\]
We use equation (5.9) and Proposition 4.1 to simplify the left side and equations (5.2) and (5.4) to simplify the right side mod $\theta$.  This yields
\[\hat{\omega}_{\#}^{\phantom{2}a}\wedge\hat{\omega}_{a}^{\phantom{2}j} = d\omega_{\#\phantom{2}\bar{\nu}}^{\phantom{2}j}\wedge\theta^{\bar{\nu}} + \hat{\omega}_{\#\phantom{2}\bar{\nu}}^{\phantom{2}j}d\theta^{\bar{\nu}} + (\hat{\omega}_{\#\phantom{2}0}^{\phantom{2}j}g_{\mu\bar{\nu}} + D_{\#}^{\phantom{2}j}g_{\mu\bar{\nu}})\theta^{\mu}\wedge\theta^{\bar{\nu}}.\]
We now only consider terms involving $\theta^{\mu}\wedge\theta^{\bar{\nu}}.$  Hence we now have
\[(\hat{\omega}_{\#\phantom{2}\mu}^{\phantom{2}a}\hat{\omega}_{a\phantom{2}\bar{\nu}}^{\phantom{2}j} - \hat{\omega}_{a\phantom{2}\mu}^{\phantom{2}j}\hat{\omega}_{\#\phantom{2}\bar{\nu}}^{\phantom{2}a})\theta^{\mu}\wedge\theta^{\bar{\nu}} = d\omega_{\#\phantom{2}\bar{\nu}}^{\phantom{2}j}\wedge\theta^{\bar{\nu}} + \hat{\omega}_{\#\phantom{2}\bar{\nu}}^{\phantom{2}j}d\theta^{\bar{\nu}} + (\hat{\omega}_{\#\phantom{2}0}^{\phantom{2}j}g_{\mu\bar{\nu}} + D_{\#}^{\phantom{2}j}g_{\mu\bar{\nu}})\theta^{\mu}\wedge\theta^{\bar{\nu}}.\]
After using the structure equation from (4.1) for $d\theta^{\bar{\nu}}$, Proposition 4.1, and simplifying, we note that $d\theta^{\bar{\alpha}} \equiv -\omega_{\bar{\nu}}^{\:\:\:\bar{\alpha}}\wedge\theta^{\bar{\nu}}$ mod $\theta$, so the coefficient of $\theta^{\mu}\wedge\theta^{\bar{\nu}}$ in the expression $\hat{\omega}_{\#\phantom{2}\bar{\alpha}}^{\phantom{2}j}d\theta^{\bar{\alpha}}$ is $-\hat{\omega}_{\#\phantom{2}\bar{\alpha}}^{\phantom{2}j}\hat{\omega}_{\bar{\nu}\phantom{2}\mu}^{\phantom{2}\bar{\alpha}}.$  Hence we are left with the equality
\[
(d\omega_{\#\phantom{2}\bar{\nu}}^{\phantom{2}j})_{\mu} - \hat{\omega}_{a\phantom{2}\bar{\nu}}^{\phantom{2}j}\hat{\omega}_{\#\phantom{2}\mu}^{\phantom{2}a} +
\hat{\omega}_{\#\phantom{2}\bar{\nu}}^{\phantom{2}a}\hat{\omega}_{a\phantom{2}\mu}^{\phantom{2}j} - \hat{\omega}_{\#\phantom{2}\bar{\alpha}}^{\phantom{2}j}\hat{\omega}_{\bar{\nu}\phantom{2}\mu}^{\phantom{2}\bar{\alpha}} = -(\hat{\omega}_{\#\phantom{2}0}^{\phantom{2}j}g_{\mu\bar{\nu}} + D_{\#}^{\phantom{2}j}g_{\mu\bar{\nu}}).\]
However, the left hand side equals $\hat{\omega}_{\#\phantom{2}\bar{\nu};\mu}^{\phantom{2}j}$, so $\hat{\omega}_{\#\phantom{2}\bar{\nu};\mu}^{\phantom{2}j}$ is a sum of multiples of the Levi form.  We now covariantly differentiate equation (5.8) and recall that $\nabla g_{\mu\bar{\nu}} = 0$ to obtain that $\omega_{\gamma_{1}\phantom{2}\gamma_{2};\gamma_{3}\ldots\gamma_{l}}^{\phantom{2}\#}\hat{\omega}_{\#\phantom{2}\bar{\nu}}^{\phantom{2}j}$ is a sum of multiples of the Levi form, so by using Lemma 2.1 as in the derivation of (5.8), we conclude $\omega_{\gamma_{1}\phantom{2}\gamma_{2};\gamma_{3}\ldots\gamma_{l}}^{\phantom{2}\#}\hat{\omega}_{\#\phantom{2}\bar{\nu}}^{\phantom{2}j}= 0$.  This proves that the first expression in (5.11) vanishes.

Now we examine the same identity but this time look at coefficients of $\theta^{\mu}\wedge{\theta}$ so we work modulo $\theta \wedge\theta^{\bar{\beta}}$ and $\theta^{\alpha}\wedge\theta^{\bar{\beta}}$.  Since $\hat{\phi}_{\#}^{\phantom{2}j} = \omega_{\#}^{\phantom{2}j} + D_{\#}^{\phantom{2}j}\theta$, we have $d\hat{\phi}_{\#}^{\phantom{2}j} \equiv d\omega_{\#}^{\phantom{2}j} + dD_{\#}^{\phantom{2}j}\wedge\theta.$  On the other hand, we use the structure equations (4.1) and simplify, yielding the identity
\[\hat{\phi}_{\#}^{\phantom{2}a}\wedge\hat{\phi}_{a}^{\phantom{2}j} \equiv d(\hat{\omega}_{\#\phantom{2}0}^{\phantom{2}j} + D_{\#}^{\phantom{2}j})\wedge\theta,\]
so we rewrite the left hand side using Proposition 4.1, simplify, and collect coefficients of $\theta^{\mu}\wedge\theta$.  This gives
\[\partial_{\mu}(\hat{\omega}_{\#\phantom{2}0}^{\phantom{2}j} + D_{\#}^{\phantom{2}j}) + \omega_{a\phantom{2}\mu}^{\phantom{2}j}(\hat{\omega}_{\#\phantom{2}0}^{\phantom{2}a} + D_{\#}^{\phantom{2}a}) -
\omega_{\#\phantom{2}\mu}^{\phantom{2}a}(\hat{\omega}_{a\phantom{2}0}^{\phantom{2}j} + D_{a}^{\phantom{2}j}) = 0\]
which implies that $(\hat{\omega}_{\#\phantom{2}0}^{\phantom{2}j} + D_{\#}^{\phantom{2}j})_{;\mu}$ is zero.  Therefore all higher order covariant derivatives in the directions $\theta^{\alpha}, \theta^{\bar{\beta}}$ are zero, so by Lemma 2.1, this implies that the second expression in equation (5.11) vanishes.  Hence we now have that $\hat{\phi}_{\#}^{\phantom{2}j} = 0.$

Since $\phi_{\alpha}^{\phantom{2}j} = 0$, we examine $d\phi_{\alpha}^{\phantom{2}j}$ and use the structure equation and our previous results to obtain
\begin{align*} 0 &= \phi_{\alpha}^{\phantom{2}A}\wedge\phi_{A}^{\phantom{2}j} + i \theta_{\alpha}\wedge\phi^{j} \nonumber \\
&= i\theta_{\alpha}\wedge(\hat{E}^{j}\theta). \nonumber
\end{align*}
This implies that $\hat{E}^{j} = 0$, so $\hat{\phi}^{j} = 0$ also.

So far we have shown that $\hat{\phi}_{\alpha}^{\phantom{2}j} = \hat{\phi}_{\#}^{\phantom{2}j} = \hat{\phi}^{j} = 0.$  We choose an adapted $Q$-frame $(Z_{\Lambda})$ on $Q$ near $f(p).$  We can choose $(Z_{\Lambda})$ corresponding to our coframe $(\theta, \theta^{A})$ such that the following relations are satisfied (see the second row of equation (4.7)).
\[\Pi_{A}^{\phantom{2}0} = -i\hat{\phi}_{A}, \qquad \Pi_{A}^{\phantom{2}B} = \hat{\phi}_{A}^{\phantom{2}B} + \delta_{A}^{B}\Pi_{0}^{\phantom{2}0}, \qquad \Pi_{A}^{\phantom{2}\hat{n}+1} = 2i\theta_{A}. \]

First, note that
\[\Pi_{j}^{\phantom{2}\hat{n}+1} = 2i\theta_{j} = 2i\theta^{\bar{A}}g_{\bar{A}j}=0\]
because $\theta^{\bar{a}} = 0$ on $M$ and $g_{\bar{\alpha}j}=0.$

Next, we see that
\[\Pi_{j}^{\phantom{2}0} = -i\hat{\phi}_{j} = -i\hat{\phi}^{\bar{A}}g_{\bar{A}j} = -i\hat{\phi}^{\bar{\alpha}}g_{\bar{\alpha}j} -i\hat{\phi}^{\bar{\#}}g_{\bar{\#}j} -i\hat{\phi}^{\bar{i}}g_{\bar{i}j}.\]
The first term in the above sum is zero because of the indices of the Levi form.  The second term is zero again because of the indices of the Levi form, due to our change of basis at the beginning of the proof.  The third term is zero because $\hat{\phi}^{\bar{i}}=0.$

Now we analyze $\Pi_{j}^{\phantom{2}\alpha}$, noting that $\delta_{j}^{\phantom{2}\alpha} = 0$ and using the symmetry relation $\hat{\phi}_{j\bar{\beta}} = -\hat{\phi}_{\bar{\beta}j}.$  We have
\[ \Pi_{j}^{\phantom{2}\alpha} = \hat{\phi}_{j}^{\phantom{2}\alpha} = \hat{\phi}_{j\bar{\beta}}g^{\bar{\beta}\alpha} = -\hat{\phi}_{\bar{\beta}j}g^{\bar{\beta}\alpha} = -\hat{\phi}_{\bar{\beta}}^{\phantom{2}\bar{A}}g_{\bar{A}j}g^{\bar{\beta}\alpha} = 0 \]
because $g_{\bar{A}j} = 0$ unless $A$ is in the range of $j$, and then $\hat{\phi}_{\bar{\alpha}}^{\phantom{2}\bar{j}} = 0.$

We perform a similar analysis of $\Pi_{j}^{\phantom{2}\#}.$
\[\Pi_{j}^{\phantom{2}\#} = \hat{\phi}_{j}^{\#} = -\hat{\phi}_{\bar{A}}^{\phantom{2}\bar{B}}g_{\bar{B}j}g^{\bar{A}\#} = 0,\]
because $g_{\bar{B}j} = 0$ unless $B$ is in the range of $j$, $g^{\bar{A}\#} = 0$ unless $A$ is in the range of $\#$, and if both of these cases occur, then $\hat{\phi}_{\bar{\#}}^{\phantom{2}\bar{j}} = 0.$

This shows that $\Pi_{j}^{\phantom{2}\Omega} = 0$ unless $\Omega \in \{n+r+1, \ldots, N\}.$  Therefore, since the Maurer-Cartan forms are defined by $dZ_{\Lambda} = \Pi_{\Lambda}^{\phantom{2}\Omega}Z_{\Omega},$ we have
\begin{equation} dZ_{i} = \Pi_{i}^{\phantom{2}j}Z_{j},
\end{equation}
expressing that the derivatives of the vectors $Z_{i}$ are linear combinations of $Z_{j}$ at each point.  The proof now concludes exactly as in \cite{EHZ}, section 9.

\section{Dimensions of $E_{k}$ for Embeddings}

We now state a theorem which relates the dimensions of the $E_{k}$ for two embeddings.  To simplify notation, we write $\omega_{\alpha}^{a}$ where $a \in \{1, \ldots, N-n\}$ rather than $\omega_{\alpha}^{a+n}$ for the second fundamental forms of the mappings.  The proof is given in section 7.

\newtheorem{Theorem 2}{Theorem}[section]
\begin{Theorem 2}
Let $M \subseteq \mathbb{C}^{n+1}$ be a smooth Levi-nondegenerate hypersurface of signature $l \le n/2$ and $p \in M$.  Let $f_{0}:M \rightarrow Q_{l}^{N_{0}}$ and $f:M \rightarrow Q_{l'}^{N}$ be smooth CR mappings that are CR transversal to $Q_{l}^{N_{0}}$ at $f_{0}(p)$ and $Q_{l'}^{N}$ at $f(p)$, respectively, and $N_{0} \le N$.    Fix an admissible coframe $(\theta, \theta^{\alpha})$ on $M$ and choose corresponding coframes $(\mathring{\theta}, \mathring{\theta}^{A})_{A=1,\ldots,N_{0}}$ and $(\hat{\theta}, \hat{\theta}^{A})_{A=1,\ldots,N}$ on $Q_{l}^{N_{0}}$ and $Q_{l'}^{N}$ adapted to $f_{0}(M)$ and $f(M)$, respectively.  Let $(\mathring{\omega}_{\gamma_{1}\:\:\gamma_{2}}^{\:\:a})_{a=1,\ldots, N_{0} - n}$ and $(\omega_{\gamma_{1}\:\:\gamma_{2}}^{\:\:a})_{a=1,\ldots,N-n}$ denote the second fundamental forms of $f_{0}$ and $f$, respectively, relative to these coframes.  Let $k \ge 2$ be an integer and assume that the spaces $\mathring{E}_{j}(q)$ and $E_{j}(q)$ for $2 \le j \le k$, are of constant dimension for $q$ near $p$.  Then for each $k$, 

\begin{description}

\item[(a)] If $l = n/2$ or $l' = N/2$ or $f$ is side preserving, and if either $(N_{0} - n) + (l' - l) < l$ or $(N - l') - (n - l) < l$, we have $dim(E_{k}) \le dim(\mathring{E}_{k}) + min(l' - l, (N-l') - (n-l)).$

\item[(b)] If $f$ is side reversing and if $l' < n$, we have $dim(E_{k}) \le dim(\mathring{E}_{k}) + l' - (n- l)$.

\end{description}

\end{Theorem 2}
We may now prove Theorem 1.2.  We use the notation of Theorem 6.1.  

\theoremstyle{definition}
\newtheorem*{remark 3}{Proof of Theorem 1.2}
\begin{remark 3}

If $l = n/2$, or $f$ is side preserving, we notice that $l' \geq l$ and $N - l' \geq n -l$ by Proposition 3.1.  Next, we apply Theorem 6.1.  Since dim$E_{k} \leq (N_{0} - n) + min(l' - l, (N-n) - (l' -l))$ for all $k$, the degeneracy of $f$ is at least $(N - n) - (N_{0} - n) - min(l' - l, (N-n) - (l' -l))$, so if $s$ denotes the degeneracy, we have $s \geq (N - N_{0}) - min(l' - l, (N-n) - (l'-l))$.  Since $(N - n) - s \leq (N_{0} - n) + (l' - l) < n$, we may apply Theorem 1 to obtain the desired result.

If $f$ is side reversing, we notice that $N - l' \geq l$ and $l' \geq n - l$ by Proposition 3.1.  We apply Theorem 6.1 again to see that the degeneracy of $f$ is at least $(N - n) - (N_{0} - n) - (l' - (n-l))$.  Denoting the degeneracy by $s$ again, we have $s \geq (N - l' - l) + (n - N_{0})$.  Since $(N - n) - s \leq (N_{0} -n) + l' - (n-l) < n$, we may apply Theorem 1 to obtain the desired result.

\end{remark 3}

A key ingredient in the proof of Theorem 6.1  is the Gauss equation for the second fundamental form of the embedding.  A more general and precise version is stated and proved in \cite{EHZ} where it appears as Theorem 2.3.  The statement here is the same as Lemma 4.3 in \cite{BEH1}.

\theoremstyle{plain}
\newtheorem{Gauss equation}[Theorem 2]{Lemma}
\begin{Gauss equation}
Let $M \subset \mathbb{C}^{n+1}$ be a smooth Levi-nondegenerate hypersurface of signature $l\le\frac{n}{2}$, $f: M \longrightarrow Q_{l'}^{N} \subset \mathbb{C}^{N+1}$ a smooth CR mapping that is CR transversal to $Q_{l'}^{N}$ along $M$, $l \leq l'$, and $\omega_{\alpha\textrm{ }\beta}^{\phantom{1}a}$ its second fundamental form.  Then,
\begin{equation*}
0 = S_{\alpha\bar{\beta}\mu\bar{\nu}} + g_{a\bar{b}}\omega_{\alpha\:\:\:\mu}^{\:\:\:a}\omega_{\bar{\beta}\:\:\:\bar{\nu}}^{\:\:\:\bar{b}} + T_{\alpha\bar{\beta}\mu\bar{\nu}},
\end{equation*}
where $S_{\alpha\bar{\beta}\mu\bar{\nu}}$ is the Chern-Moser pseudoconformal curvature of $M$ and $T_{\alpha\bar{\beta}\mu\bar{\nu}}$ is a conformally flat tensor.
\end{Gauss equation}

We shall need the following lemma regarding conformal flatness of certain covariant derivatives of the second fundamental form.  This lemma appears with proof as Lemma 4.1 in \cite{BEH1}.

\newtheorem{Lemma 1}[Theorem 2]{Lemma}
\begin{Lemma 1}
Let $M$, $f$, and $\omega_{\alpha\textrm{ }\beta}^{\phantom{1}a}$ be as in Lemma 6.2.  Then the covariant derivative tensor $\omega_{\alpha\:\:\:\beta;\bar{\gamma}}^{\:\:a}$ is conformally flat.
\end{Lemma 1}

It will also be necessary to know how covariant derivatives of the second fundamental form commute.  Given a CR embedding $f:M \rightarrow \hat{M}$, we now recall some facts about the pseudoconformal connection on $\hat{M}$ pulled back to $M$.  Suppose $(\theta, \theta^{A})$ is an adapted coframe for the pair $(M, \hat{M})$.  We use the same notation as in the Preliminaries section.  We denote with a $\hat{}$ the pseudoconformal connection forms on $\hat{M}$ pulled back to $M$, where the indices run from $1$ to $\hat{n}$.  Recall that $(\omega, \omega^{\alpha}, \omega^{\bar{\alpha}}) = (\hat{\omega}, \hat{\omega}^{\alpha}, \hat{\omega}^{\bar{\alpha}}) = (\theta, \theta^{\alpha}, \theta^{\bar{\alpha}})$ and $\hat{\omega}^{a} = 0$ on $M$.  We do not expect $(\phi_{\beta}^{\:\:\:\alpha}, \phi^{\alpha}, \psi)$ and $(\hat{\phi}_{\beta}^{\:\:\:\alpha}, \hat{\phi}^{\alpha}, \hat{\psi})$ to be equal, but since $\hat{\omega}_{\beta}^{\:\:\:\alpha} = \omega_{\beta}^{\:\:\:\alpha}$ and $\hat{\tau}^{\alpha} = \tau^{\alpha}$, Proposition 4.1 implies
\begin{equation}
\hat{\phi}_{\beta}^{\:\:\:\alpha} = \phi_{\beta}^{\:\:\:\alpha} + C_{\beta}^{\:\:\:\alpha}\theta, \qquad  \hat{\phi}^{\alpha} =  \phi^{\alpha} + C_{\mu}^{\:\:\:\alpha}\theta^{\mu} + F^{\alpha}\theta, \qquad \hat{\psi} = \psi + iF_{\mu}\theta^{\mu} - iF_{\bar{\nu}}\theta^{\bar{\nu}} + A\theta
\end{equation}
where
\begin{equation*}
C_{\beta}^{\:\:\:\alpha} := \hat{D}_{\beta}^{\:\:\:\alpha} - D_{\beta}^{\:\:\:\alpha}, \qquad F^{\alpha} := \hat{E}^{\alpha} - E^{\alpha}, \qquad A:= \hat{B} - B
\end{equation*}
and $\hat{D}_{\beta}^{\:\:\:\alpha}, \hat{E}^{\alpha}, \hat{B}$ are the analogues for $\hat{M}$ of the functions from Proposition 4.1 restricted to $M$.  We also record the following expression for $C_{\alpha\bar{\beta}}$ which appears as equation (6.8) in \cite{EHZ}.
\begin{equation}
C_{\alpha\bar{\beta}} = \frac{i(\hat{S}_{a\:\:\:\alpha\bar{\beta}}^{\:\:\:a} + \omega_{\mu\:\:\:\alpha}^{\:\:\:a}\omega_{\:\:\:a\bar{\beta}}^{\mu})}{n+2} - \frac{i(\hat{S}_{a\:\:\:\mu}^{\:\:\:a \:\:\: \mu} + \omega_{\mu\:\:\:\nu}^{\:\:\:a}\omega_{\:\:\:a}^{\mu\:\:\:\nu})g_{\alpha\bar{\beta}}}{2(n+1)(n+2)}. 
\end{equation}

The following is a more specific version of Lemma 4.2 in \cite{BEH1}, where we give an explicit formula for the part which is not conformally flat.
\newtheorem{Lemma 2}[Theorem 2]{Lemma}
\begin{Lemma 2}
Let $M$, $f$, and $\omega_{\alpha\textrm{ }\beta}^{\phantom{1}a}$ be as in Lemma 6.2, and $p\in M$.  Then for any $s \ge 2$, we have
\begin{equation}
\omega_{\gamma_{1} \phantom{1} \gamma_{2};\gamma_{3}\ldots\gamma_{s}\alpha\bar{\beta}}^{\phantom{1}a} - \omega_{\gamma_{1} \phantom{1} \gamma_{2};\gamma_{3}\ldots\gamma_{s}\bar{\beta}\alpha}^{\phantom{1}a} \equiv \sum_{j=1}^{s}d_{\alpha\bar{\beta}}(\omega_{\gamma_{j}}^{\phantom{1}\mu})\omega_{\gamma_{1} \phantom{1} \gamma_{2};\gamma_{3}\ldots\gamma_{j-1}\mu\gamma_{j+1}\ldots\gamma_{s}}^{\phantom{1}a} - C_{\phantom{1}\alpha\bar{\beta}c}^{a}\omega_{\gamma_{1} \phantom{1} \gamma_{2};\gamma_{3}\ldots\gamma_{s}}^{\phantom{1}c}
\end{equation}
where equivalence is modulo a conformally flat tensor, $d_{\alpha\bar{\beta}}(\omega_{\gamma_{j}}^{\phantom{1}\mu})$ is the coefficient of $\theta^{\alpha}\wedge\theta^{\bar{\beta}}$ in $d\omega_{\gamma_{j}}^{\phantom{1}\mu}$, and $C_{\phantom{1}\alpha\bar{\beta}c}^{a}$ is given by
\begin{equation*}
C_{\phantom{1}\alpha\bar{\beta}c}^{a} \equiv \omega_{\phantom{15}\alpha}^{\bar{\rho}a}\omega_{\bar{\rho}c\bar{\beta}} + i\delta_{c}^{a}\hat{D}_{\bar{\beta}\alpha}.
\end{equation*}

\begin{proof}
We use the pseudoconformal connections introduced in section 4.  We observe that the left hand side of (6.3) is a tensor, hence it is enough to show (6.3) at each fixed $p \in M$ with respect to any choice of adapted coframe near $p$.  By making a unitary change of coframe $\theta^{\alpha} \rightarrow u_{\beta}^{\phantom{2}\alpha}\theta^{\beta}$ and $\theta^{a} \rightarrow u_{b}^{\phantom{2}a}\theta^{b}$ in the tangential and normal directions, we may choose an adapted coframe near $p$ such that $\omega_{\alpha}^{\phantom{1}\beta}(p) = \omega_{a}^{\phantom{1}b}(p) = 0$ (c.f. Lemma 2.1 in \cite{Le88}.
In this coframe, the left hand side at $p$ is equivalent, modulo a conformally flat tensor, to the coefficient in front of $\theta^{\alpha}\wedge\theta^{\bar{\beta}}$ in the expression
\begin{equation*}
\sum_{j=1}^{s}d\omega_{\gamma_{j}}^{\phantom{1}\mu}\omega_{\gamma_{1} \phantom{1} \gamma_{2};\gamma_{3}\ldots\gamma_{j-1}\mu\gamma_{j+1}\ldots\gamma_{s}}^{\phantom{1}a} -
\omega_{\gamma_{1} \phantom{1} \gamma_{2};\gamma_{3}\ldots\gamma_{s}}^{\phantom{1}c}d\hat{\phi}_{c}^{\phantom{1}a}.
\end{equation*}
Hence we would like to show that the coefficient in front of $\theta^{\alpha}\wedge\theta^{\bar{\beta}}$ in $d\hat{\phi}_{c}^{\phantom{1}a}$ has the form of the $C_{\phantom{1}\alpha\bar{\beta}c}^{a}$ given in the statement of the Lemma.

Note that we may work mod $\theta$ because we are only looking for  the coefficient in front of $\theta^{\alpha}\wedge\theta^{\bar{\beta}}$.  The structure equations (4.1) give
\begin{equation*}
d\hat{\phi}_{c}^{\phantom{1}a} \equiv \hat{\phi}_{c}^{\phantom{1}\rho}\wedge \hat{\phi}_{\rho}^{\phantom{1}a} - i\delta_{c}^{\phantom{1}a}\hat{\phi}_{\mu}\wedge \theta^{\mu}.
\end{equation*}
Now using Proposition 4.1, and equation (3.3), $\omega_{B\bar{A}} + \omega_{\bar{A}B} = 0$, we have
\begin{eqnarray}
 \hat{\phi}_{c}^{\phantom{1}\rho} \equiv \omega_{c}^{\phantom{1}\rho} &\equiv& \omega_{c\bar{A}}g^{\rho\bar{A}} \nonumber \\
& \equiv &  - \omega_{\bar{A}c}g^{\rho\bar{A}}  \nonumber \\
& \equiv & - \overline{\omega_{A}^{\phantom{1}B} g_{B\bar{c}}g^{A\bar{\rho}}} \nonumber \\
& \equiv & - \overline{\hat{\phi}_{\rho}^{\phantom{1}d}}g^{\bar{\rho}\mu}g_{c\bar{d}}. \nonumber
\end{eqnarray}
Hence the coefficient of $\theta^{\alpha}\wedge\theta^{\bar{\beta}}$ from $\hat{\phi}_{c}^{\phantom{1}\rho}\wedge \hat{\phi}_{\rho}^{\phantom{1}a}$ is $\omega_{\phantom{15}\alpha}^{\bar{\rho}a}\omega_{\bar{\rho}c\bar{\beta}}.$

We next examine $i\delta_{c}^{\phantom{1}a}\hat{\phi}_{\mu}\wedge \theta^{\mu}$ and work mod $\theta$.  We notice that
\begin{equation*}
i\delta_{c}^{\phantom{1}a}\hat{\phi}_{\mu}\wedge \theta^{\mu} = i\delta_{c}^{\phantom{1}a}\theta^{\mu}\wedge(g_{\mu\bar{A}}\hat{\phi}^{\bar{A}}) = i\delta_{c}^{\phantom{1}a}\theta^{\mu}\wedge(g_{\mu\bar{\sigma}}\hat{\phi}^{\bar{\sigma}})
\end{equation*}
due to the form of the matrix $(g)$.  We substitute for $\hat{\phi}^{\bar{\sigma}}$ using equation (6.1) and use Proposition 4.1 to obtain
\begin{equation*}
i\delta_{c}^{\phantom{1}a}\theta^{\mu}\wedge(g_{\mu\bar{\sigma}}(\phi^{\bar{\sigma}} + C_{\bar{\rho}}^{\phantom{1}\bar{\sigma}}\theta^{\bar{\rho}} + F^{\bar{\sigma}}\theta))
\equiv
i\delta_{c}^{\phantom{1}a}\theta^{\mu}\wedge(g_{\mu\bar{\sigma}}(\tau^{\bar{\sigma}} + (D_{\bar{\rho}}^{\phantom{1}\bar{\sigma}} + C_{\bar{\rho}}^{\phantom{1}\bar{\sigma}})\theta^{\bar{\rho}}).
\end{equation*}
We notice that by equation (3.2), $\tau^{\bar{\sigma}}$ will be a combination of only forms like $\theta^{\gamma}$, so we may ignore it when searching for coefficients of $\theta^{\alpha}\wedge\theta^{\bar{\beta}}$.  We also note that $\hat{D}_{\bar{\beta}}^{\phantom{1}\bar{\alpha}} = C_{\bar{\beta}}^{\phantom{1}\bar{\alpha}} + D_{\bar{\beta}}^{\phantom{1}\bar{\alpha}}$ by (6.1), so after lowering an index, we find the coefficient of $i\delta_{c}^{\phantom{1}a}\hat{\phi}_{\mu}\wedge \theta^{\mu}$ in front of $\theta^{\alpha}\wedge\theta^{\bar{\beta}}$ is exactly $(\omega_{\phantom{15}\alpha}^{\bar{\rho}a}\omega_{\bar{\rho}c\bar{\beta}} + i\delta_{c}^{a}\hat{D}_{\bar{\beta}\alpha})$, as desired.

\end{proof}
\end{Lemma 2}

The following linear algebra Lemma will be useful

\newtheorem{Linear Algebra Lemma}[Theorem 2]{Lemma}
\begin{Linear Algebra Lemma}

Let $\{w_{1}, \ldots, w_{m}\}$ and $\{v_{1}, \ldots, v_{m}\}$ be vectors in $\mathbb{C}^{n}$ such that $\langle w_{i}, w_{j} \rangle' = \langle v_{i}, v_{j} \rangle.$ Here    $\langle x, y \rangle$ denotes the standard inner product on $\mathbb{C}^{n}$ and $\langle x, y \rangle' = y^{*}I_{k}x$, where $I_{k}$ is the $n\times n$ diagonal matrix with first $k$ entries equal to $-1$ and remaining $n-k$ entries equal to $1$.  Let $W = \textrm{span}\{w_{1}, \ldots, w_{m}\}$ and $V = \textrm{span}\{v_{1}, \ldots, v_{m}\}$, then $dim(W) \le dim(V) + min(k, n-k)$.

\begin{proof}
Let $w_{i_{1}}, \ldots, w_{i_{k}}$ be a basis for $W$ and define a linear map $\phi$ from $W$ to $V$ by $\phi(w_{i_{j}}) = v_{i_{j}}$.  Suppose $x = \sum_{r=1}^{k} a_{r}w_{i_{r}}$ is in the kernel of $\phi$, so $\sum_{r=1}^{k} a_{r}v_{i_{r}} = 0$.  If we let $a$ denote the coordinate vector of $x$, we have
\begin{math}
\langle x, x \rangle' = a^{*}\big(\langle w_{i_{r}}, w_{i_{s}} \rangle' \big)a = a^{*}\big(\langle v_{i_{r}}, v_{i_{s}} \rangle \big)a = \langle \phi(x), \phi(x) \rangle = 0.
\end{math}
This shows that the kernel of $\phi$ is an isotropic subspace of $\mathbb{C}^{n}$ with respect to $\langle \cdot, \cdot \rangle'$, which implies that the dimension of the kernel is at most $min(k, n-k)$.  The result follows by standard linear algebra.

\end{proof}
\end{Linear Algebra Lemma}

\section{Proof of Theorem 6.1}

We first prove by induction that  for all $j, k \ge 2$, we have
\begin{equation}
g_{a\bar{b}}\omega_{\gamma_{1} \: \gamma_{2} ; \gamma_{3} \ldots \gamma_{j}}^{\:\:\:a}\omega_{\bar{\alpha}_{1}\:\bar{\alpha}_{2};\bar{\alpha}_{3}\ldots\bar{\alpha}_{k}}^{\:\:\:\bar{b}} \equiv \:\: \mathring{g}_{a\bar{b}}\mathring{\omega}_{\gamma_{1} \: \gamma_{2} ; \gamma_{3} \ldots \gamma_{j}}^{\:\:\:a}\mathring{\omega}_{\bar{\alpha}_{1}\:\bar{\alpha}_{2};\bar{\alpha}_{3}\ldots\bar{\alpha}_{k}}^{\:\:\:\bar{b}},
\end{equation}
where equivalence here and in the rest of the proof means that the sides of the equation differ by a conformally flat tensor.  We then show that such conformal equivalence is in fact equality and apply Lemma 6.5.  We induct on the sum of the indices.  By subtracting the Gauss equation for $\omega_{\gamma_{1}\:\gamma_{2}}^{\:\:\:a}$ from the corresponding one for $\mathring{\omega}_{\gamma_{1}\:\gamma_{2}}^{\:\:\:a}$, we obtain
\begin{equation}
g_{a\bar{b}}\omega_{\gamma_{1} \: \gamma_{2}}^{\:\:\:a}\omega_{\bar{\alpha}_{1}\: \bar{\alpha}_{2}}^{\:\:\:\bar{b}} \equiv \mathring{g}_{a\bar{b}}\mathring{\omega}_{\gamma_{1} \: \gamma_{2}}^{\:\:\:a}\mathring{\omega}_{\bar{\alpha}_{1}\: \bar{\alpha}_{2}}^{\:\:\:\bar{b}},
\end{equation}
since the pseudoconformal curvature tensor $S_{\gamma_{1}\bar{\alpha}_{1}\gamma_{2}\bar{\alpha}_{2}}$ is computed using the same coframe $(\theta, \theta^{\alpha})$.  This establishes the base step of the induction.  We now assume equation (7.1) with $j + k \leq p$, and we wish to show the same where $j + k = p + 1$.  We will demonstrate the case where $k$ increases by $1$.  The case where $j$ increases is similar and left to the reader.  We differentiate both sides of (7.1) in the $\theta^{\bar{\gamma}_{k+1}}$ direction, note that covariant derivatives of conformally flat tensors are conformally flat, and obtain
\begin{multline*}
g_{a\bar{b}}\omega_{\gamma_{1} \: \gamma_{2} ; \gamma_{3} \ldots \gamma_{j}\bar{\alpha}_{k+1}}^{\:\:\:a}\omega_{\bar{\alpha}_{1}\:\bar{\alpha}_{2};\bar{\alpha}_{3}\ldots\bar{\alpha}_{k}}^{\:\:\:\bar{b}} +
g_{a\bar{b}}\omega_{\gamma_{1} \: \gamma_{2} ; \gamma_{3} \ldots \gamma_{j}}^{\:\:\:a}\omega_{\bar{\alpha}_{1}\:\bar{\alpha}_{2};\bar{\alpha}_{3}\ldots\bar{\alpha}_{k}\bar{\alpha}_{k+1}}^{\:\:\:\bar{b}} \equiv
\\
\mathring{g}_{a\bar{b}}\mathring{\omega}_{\gamma_{1} \: \gamma_{2} ; \gamma_{3} \ldots \gamma_{j}\bar{\alpha}_{k+1}}^{\:\:\:a}\mathring{\omega}_{\bar{\alpha}_{1}\:\bar{\alpha}_{2};\bar{\alpha}_{3}\ldots\bar{\alpha}_{k}}^{\:\:\:\bar{b}} +
\mathring{g}_{a\bar{b}}\mathring{\omega}_{\gamma_{1} \: \gamma_{2} ; \gamma_{3} \ldots \gamma_{j}}^{\:\:\:a}\mathring{\omega}_{\bar{\alpha}_{1}\:\bar{\alpha}_{2};\bar{\alpha}_{3}\ldots\bar{\alpha}_{k}\bar{\alpha}_{k+1}}^{\:\:\:\bar{b}}.
\end{multline*}
The next lemma shows the equivalence of the first terms on each side of the above equation.  We then subtract to finish the induction and demonstrate equation (7.1) for all $j,k \geq 2$.

\newtheorem{Lemma 3}{Lemma}[section]
\begin{Lemma 3}
With the same setup as above, we have
\begin{equation*}
g_{a\bar{b}}\omega_{\gamma_{1} \: \gamma_{2} ; \gamma_{3} \ldots \gamma_{j}\bar{\alpha}_{k+1}}^{\:\:\:a}\omega_{\bar{\alpha}_{1}\:\bar{\alpha}_{2};\bar{\alpha}_{3}\ldots\bar{\alpha}_{k}}^{\:\:\:\bar{b}} \equiv \:\:
\mathring{g}_{a\bar{b}}\mathring{\omega}_{\gamma_{1} \: \gamma_{2} ; \gamma_{3} \ldots \gamma_{j}\bar{\alpha}_{k+1}}^{\:\:\:a}\mathring{\omega}_{\bar{\alpha}_{1}\:\bar{\alpha}_{2};\bar{\alpha}_{3}\ldots\bar{\alpha}_{k}}^{\:\:\:\bar{b}}
\end{equation*}
\begin{proof}
We induct on $s$, the position of the index $\bar{\alpha}_{k+1}$.  When $s = 3$, we need
\begin{equation*}
g_{a\bar{b}}\omega_{\gamma_{1} \: \gamma_{2} ; \bar{\alpha}_{k+1}\gamma_{3}\ldots\gamma_{j}}^{\:\:\:a}\omega_{\bar{\alpha}_{1}\:\bar{\alpha}_{2};\bar{\alpha}_{3}\ldots\bar{\alpha}_{k}}^{\:\:\:\bar{b}} \equiv \:\:
\mathring{g}_{a\bar{b}}\mathring{\omega}_{\gamma_{1} \: \gamma_{2} ; \bar{\alpha}_{k+1}\gamma_{3}\ldots\gamma_{j}}^{\:\:\:a}\mathring{\omega}_{\bar{\alpha}_{1}\:\bar{\alpha}_{2};\bar{\alpha}_{3}\ldots\bar{\alpha}_{k}}^{\:\:\:\bar{b}}.
\end{equation*}
This follows immediately from Lemma 6.3, which implies that both sides are conformally flat.

We assume that the desired equivalence holds for $s = r$, where $r \leq j$, that is,
\begin{equation}
g_{a\bar{b}}\omega_{\gamma_{1} \: \gamma_{2} ; \gamma_{3} \ldots \gamma_{r-1}\bar{\alpha}_{k+1}\gamma_{r}\ldots\gamma_{j}}^{\:\:\:a}\omega_{\bar{\alpha}_{1}\:\bar{\alpha}_{2};\bar{\alpha}_{3}\ldots\bar{\alpha}_{k}}^{\:\:\:\bar{b}} \equiv \:\:
\mathring{g}_{a\bar{b}}\mathring{\omega}_{\gamma_{1} \: \gamma_{2} ; \gamma_{3} \ldots \gamma_{r-1}\bar{\alpha}_{k+1}\gamma_{r}\ldots\gamma_{j}}^{\:\:\:a}\mathring{\omega}_{\bar{\alpha}_{1}\:\bar{\alpha}_{2};\bar{\alpha}_{3}\ldots\bar{\alpha}_{k}}^{\:\:\:\bar{b}},
\end{equation}
and we would like to show the same when $s = r + 1$:
\begin{equation}
g_{a\bar{b}}\omega_{\gamma_{1} \: \gamma_{2} ; \gamma_{3} \ldots \gamma_{r}\bar{\alpha}_{k+1}\gamma_{r+1}\ldots\gamma_{j}}^{\:\:\:a}\omega_{\bar{\alpha}_{1}\:\bar{\alpha}_{2};\bar{\alpha}_{3}\ldots\bar{\alpha}_{k}}^{\:\:\:\bar{b}} \equiv \:\:
\mathring{g}_{a\bar{b}}\mathring{\omega}_{\gamma_{1} \: \gamma_{2} ; \gamma_{3} \ldots \gamma_{r}\bar{\alpha}_{k+1}\gamma_{r+1}\ldots\gamma_{j}}^{\:\:\:a}\mathring{\omega}_{\bar{\alpha}_{1}\:\bar{\alpha}_{2};\bar{\alpha}_{3}\ldots\bar{\alpha}_{k}}^{\:\:\:\bar{b}}.
\end{equation}

By Lemma 6.4, we have
\begin{multline}
\omega_{\gamma_{1} \: \gamma_{2};\gamma_{3}\ldots\gamma_{r}\bar{\alpha}_{k+1}}^{\:\:\:a} \equiv \omega_{\gamma_{1} \: \gamma_{2};\gamma_{3}\ldots\bar{\alpha}_{k+1}\gamma_{r}}^{\:\:\:a}
+ \sum_{q=1}^{r-1}d_{\gamma_{r}\bar{\alpha}_{k+1}}(\omega_{\gamma_{q}}^{\:\:\:\mu})\omega_{\gamma_{1}\:\gamma_{2};\gamma_{3}\ldots\gamma_{q-1}\mu\gamma_{q+1}\ldots\gamma_{r-1}}^{\:\:\:a} \\
- (g^{\sigma\bar{\rho}}g_{c\bar{d}}\omega_{\sigma\:\:\gamma_{r}}^{\:\:\:a}\omega_{\bar{\rho}\:\:\bar{\alpha}_{k+1}}^{\:\:\:\bar{d}})\omega_{\gamma_{1}\:\gamma_{2};\gamma_{3}\ldots\gamma_{r-1}}^{\:\:\:c}
- i(\hat{D}_{\bar{\alpha}_{k+1}\gamma_{r}})\omega_{\gamma_{1}\:\gamma_{2};\gamma_{3}\ldots\gamma_{r-1}}^{\:\:\:a}
\end{multline}
We take covariant derivatives of both sides of (7.5) in the $\theta^{\gamma_{r+1}},\ldots \theta^{\gamma_{j}}$ directions successively, multiply by $g_{a\bar{b}}\omega_{\bar{\alpha}_{1}\:\bar{\alpha}_{2};\bar{\alpha}_{3}\ldots\bar{\alpha}_{k}}^{\:\:\:\bar{b}}$, and analyze each term on the right side of the resulting equation.  We will show that by (7.3), and equation (7.1) with $j + k \leq p$ (the induction hypotheses in the proof of Lemma 7.1 and the proof of Theorem 6.1 respectively), each such term must be conformally equivalent to the corresponding term with the ring superscript.  This will demonstrate (7.4) and hence conclude the proof of Lemma 7.1.  This is  because we may also apply Lemma 6.4 to $\mathring{\omega}_{\gamma_{1} \: \gamma_{2};\gamma_{3}\ldots\gamma_{r}\bar{\alpha}_{k+1}}^{\:\:\:a}$, take covariant derivatives in the $\theta^{\gamma_{r+1}},\ldots \theta^{\gamma_{j}}$ directions, and multiply by $\mathring{g}_{a\bar{b}}\mathring{\omega}_{\bar{\alpha}_{1}\:\bar{\alpha}_{2};\bar{\alpha}_{3}\ldots\bar{\alpha}_{k}}^{\:\:\:\bar{b}}$. 

After taking covariant derivatives and multiplying by $g_{a\bar{b}}\omega_{\bar{\alpha}_{1}\:\bar{\alpha}_{2};\bar{\alpha}_{3}\ldots\bar{\alpha}_{k}}^{\:\:\:\bar{b}}$, the first term on the right side of (7.5) will be $g_{a\bar{b}}\omega_{\gamma_{1} \: \gamma_{2};\gamma_{3}\ldots\gamma_{r-1}\bar{\alpha}_{k+1}\gamma_{r}\ldots\gamma_{j}}^{\:\:\:a}\omega_{\bar{\alpha}_{1}\:\bar{\alpha}_{2};\bar{\alpha}_{3}\ldots\bar{\alpha}_{k}}^{\:\:\:\bar{b}}$ which is conformally equivalent to the same term with the ring superscript by (7.3).  We also notice that after taking covariant derivatives and multiplying by $g_{a\bar{b}}\omega_{\bar{\alpha}_{1}\:\bar{\alpha}_{2};\bar{\alpha}_{3}\ldots\bar{\alpha}_{k}}^{\:\:\:\bar{b}}$, the second term on the right side of equation (7.5) yields many terms, each of which is a product of covariant derivatives of $d_{\gamma_{r}\bar{\alpha}_{k+1}}(\omega_{\gamma_{q}}^{\:\:\:\mu})$, covariant derivatives of $\omega_{\gamma_{1}\:\gamma_{2};\gamma_{3}\ldots\gamma_{q-1}\mu\gamma_{q+1}\ldots\gamma_{r-1}}^{\:\:\:a}$, and $g_{a\bar{b}}\omega_{\bar{\alpha}_{1}\:\bar{\alpha}_{2};\bar{\alpha}_{3}\ldots\bar{\alpha}_{k}}^{\:\:\:\bar{b}}$.  We notice that expressions of the form $d_{\gamma_{r}\bar{\alpha}_{k+1}}(\omega_{\gamma_{q}}^{\:\:\:\mu})$ are intrinsic to the manifold $M$ and thus all covariant derivatives will be the same as those with the ring superscript.  Also, covariant derivatives of $\omega_{\gamma_{1}\:\gamma_{2};\gamma_{3}\ldots\gamma_{q-1}\mu\gamma_{q+1}\ldots\gamma_{r-1}}^{\:\:\:a}$ multiplied by $g_{a\bar{b}}\omega_{\bar{\alpha}_{1}\:\bar{\alpha}_{2};\bar{\alpha}_{3}\ldots\bar{\alpha}_{k}}^{\:\:\:\bar{b}}$ will be the same as those with the ring superscript by (7.1), since $r + k \leq j + k \leq p$.

The third term on the right side of equation (7.5) can be written as
\begin{equation*}
- (g^{\sigma\bar{\rho}}\omega_{\sigma\:\:\gamma_{r}}^{\:\:\:a})(g_{c\bar{d}}\omega_{\gamma_{1}\:\gamma_{2};\gamma_{3}\ldots\gamma_{r-1}}^{\:\:\:c}\omega_{\bar{\rho}\:\:\bar{\alpha}_{k+1}}^{\:\:\:\bar{d}}).
\end{equation*}
We observe that $(g_{c\bar{d}}\omega_{\gamma_{1}\:\gamma_{2};\gamma_{3}\ldots\gamma_{r-1}}^{\:\:\:c}\omega_{\bar{\rho}\:\:\bar{\alpha}_{k+1}}^{\:\:\:\bar{d}})$ is conformally equivalent to the same with the ring superscript by (7.1), and hence covariant derivatives will be also.  Also, taking covariant derivatives of the term $g^{\sigma\bar{\rho}}\omega_{\sigma\:\:\gamma_{r}}^{\:\:\:a}$ and multiplying by $g_{a\bar{b}}\omega_{\bar{\alpha}_{1}\:\bar{\alpha}_{2};\bar{\alpha}_{3}\ldots\bar{\alpha}_{k}}^{\:\:\:\bar{b}}$ will yield terms conformally equivalent to those with the ring superscript, again by (7.1) and because $g^{\sigma\bar{\rho}}$ is intrinsic to $M$.

In the last term on the right side of equation (7.5), we first show that $\hat{D}_{\bar{\alpha}_{k+1}\gamma_{r}}$ is conformally equivalent to the same with the ring superscript.  Observe that by equation (6.2), we have
\begin{equation*}
C_{\alpha\bar{\beta}} = \frac{i}{n+2}[\omega_{\mu\:\:\alpha}^{\:a}\omega_{\:\:a\bar{\beta}}^{\mu} - \frac{g_{\alpha\bar{\beta}}}{2(n+1)}\omega_{\mu\:\:\nu}^{\:a}\omega_{\:\:a}^{\mu\:\nu}].
\end{equation*}
Here we have used the vanishing of the pseudoconformal curvature of the target hyperquadric.  We may write this as
\begin{equation*}
C_{\alpha\bar{\beta}} = \frac{i}{n+2}[g^{\mu\bar{\nu}}(g_{a\bar{b}}\omega_{\mu\:\:\alpha}^{\:a}\omega_{\bar{\nu}\:\:\bar{\beta}}^{\:\bar{b}}) - \frac{g_{\alpha\bar{\beta}}g^{\mu\bar{\sigma}}g^{\nu\bar{\beta}}}{2(n+1)}(g_{a\bar{b}}\omega_{\mu\:\:\nu}^{\:a}\omega_{\bar{\sigma}\:\:\bar{\beta}}^{\:\bar{b}})].
\end{equation*}
Equation (7.2) implies conformal equivalence of both terms of the form $(g_{a\bar{b}}\omega_{\mu\:\:\alpha}^{\:a}\omega_{\bar{\nu}\:\:\bar{\beta}}^{\:\bar{b}})$ with the corresponding terms with superscripts.  Since $\hat{D}_{\bar{\beta}}^{\phantom{1}\bar{\alpha}} = C_{\bar{\beta}}^{\phantom{1}\bar{\alpha}} + D_{\bar{\beta}}^{\phantom{1}\bar{\alpha}}$ (see (6.1)), and the term $D_{\bar{\beta}}^{\phantom{1}\bar{\alpha}}$ is intrinsic to $M$, we have that $\hat{D}_{\bar{\alpha}_{k+1}\gamma_{r}}$ is conformally equivalent to its counterpart with the ring superscript.

Now we observe that after taking covariant derivatives and multiplying by $g_{a\bar{b}}\omega_{\bar{\alpha}_{1}\:\bar{\alpha}_{2};\bar{\alpha}_{3}\ldots\bar{\alpha}_{k}}^{\:\:\:\bar{b}}$ in the last term on the right side of equation (7.5), every resulting term will be a product of derivatives of $\hat{D}_{\bar{\alpha}_{k+1}\gamma_{r}}$, derivatives of $\omega_{\gamma_{1}\:\gamma_{2};\gamma_{3}\ldots\gamma_{r-1}}^{\:\:\:a}$, and $g_{a\bar{b}}\omega_{\bar{\alpha}_{1}\:\bar{\alpha}_{2};\bar{\alpha}_{3}\ldots\bar{\alpha}_{k}}^{\:\:\:\bar{b}}$.  The derivatives of $\hat{D}_{\bar{\alpha}_{k+1}\gamma_{r}}$ will be conformally equivalent to the same with the ring superscript, as explained above, and the remaining terms will be conformally equivalent to their counterparts with the ring superscript by (7.1).  This concludes the proof of Lemma 7.1

\end{proof}
\end{Lemma 3}

We now return to the proof of Theorem 6.1.  We have shown that
\begin{equation*}
g_{a\bar{b}}\omega_{\gamma_{1} \: \gamma_{2} ; \gamma_{3} \ldots \gamma_{j}}^{\:\:\:a}\omega_{\bar{\alpha}_{1}\:\bar{\alpha}_{2};\bar{\alpha}_{3}\ldots\bar{\alpha}_{k}}^{\:\:\:\bar{b}} \equiv \:\: \mathring{g}_{a\bar{b}}\mathring{\omega}_{\gamma_{1} \: \gamma_{2} ; \gamma_{3} \ldots \gamma_{j}}^{\:\:\:a}\mathring{\omega}_{\bar{\alpha}_{1}\:\bar{\alpha}_{2};\bar{\alpha}_{3}\ldots\bar{\alpha}_{k}}^{\:\:\:\bar{b}}
\end{equation*}
where the equivalence is modulo a conformally flat tensor.  Our next step is to show that this equivalence is in fact equality.  We will demonstrate this equality in the case where $l = n/2$ or $f$ is side-preserving.   To do this, we make use of Lemmas 2.1 and 2.2.  We first show equality in the case where $j = k$ using Lemma 2.2.  At the end of the proof we mention the side-reversing case.

First, suppose $(N_{0} - n) + (l' - l) < l$ and consider the following expression
\begin{equation}
g_{a\bar{b}}\omega_{\gamma_{1} \: \gamma_{2} ; \gamma_{3} \ldots \gamma_{k}}^{\:\:\:a}\omega_{\bar{\alpha}_{1}\:\bar{\alpha}_{2};\bar{\alpha}_{3}\ldots\bar{\alpha}_{k}}^{\:\:\:\bar{b}} - \mathring{g}_{a\bar{b}}\mathring{\omega}_{\gamma_{1} \: \gamma_{2} ; \gamma_{3} \ldots \gamma_{k}}^{\:\:\:a}\mathring{\omega}_{\bar{\alpha}_{1}\:\bar{\alpha}_{2};\bar{\alpha}_{3}\ldots\bar{\alpha}_{k}}^{\:\:\:\bar{b}} \equiv \:\: 0.
\end{equation}
Let $\zeta := (\zeta^{1}, \ldots, \zeta^{n})$, multiply equation (7.6) by $\zeta^{\gamma_{1}}\zeta^{\bar{\alpha}_{1}}\ldots\zeta^{\gamma_{k}}\zeta^{\bar{\alpha}_{k}}$ and sum.  Since the right side of (7.6) is conformally flat, we have
\begin{equation*}
-\sum_{a=1}^{l' - l}|\omega^{a}(\zeta)|^{2} - \sum_{b=1}^{N_{0}-n}|\mathring{\omega}^{b}(\zeta)|^{2} + \sum_{a= l' - l + 1}^{N - n}|\omega^{a}(\zeta)|^{2} = A(\zeta, \bar{\zeta})\bigg( -\sum_{i=1}^{l}|\zeta^{i}|^{2} + \sum_{j = l + 1}^{n}|\zeta^{j}|^{2}\bigg),
\end{equation*}
where $\omega^{a}(\zeta) = \omega_{\gamma_{1} \: \gamma_{2} ; \gamma_{3} \ldots \gamma_{k}}^{\:\:\:a}\zeta^{\gamma_{1}}\ldots\zeta^{\gamma_{k}}$, $\mathring{\omega}^{b}(\zeta) = \mathring{\omega}_{\gamma_{1} \: \gamma_{2} ; \gamma_{3} \ldots \gamma_{k}}^{\:\:\:b}\zeta^{\gamma_{1}}\ldots\zeta^{\gamma_{k}}$, and $A(\zeta, \bar{\zeta})$ is a polynomial is $\zeta$ and $\bar{\zeta}$.  Since we have $(N_{0} - n) + (l' - l) < l$, Lemma 2.2 implies that $A(\zeta,\bar{\zeta})$ is identically zero, so we have the desired equality, which we may rewrite as
\begin{equation}
\sum_{a= l' - l + 1}^{N - n}|\omega^{a}(\zeta)|^{2} = \sum_{a=1}^{l' - l}|\omega^{a}(\zeta)|^{2} + \sum_{b=1}^{N_{0}-n}|\mathring{\omega}^{b}(\zeta)|^{2}.
\end{equation}  
Now suppose that $(N - l') - (n - l) < l$.  We consider the expression 
\begin{equation*}
\mathring{g}_{a\bar{b}}\mathring{\omega}_{\gamma_{1} \: \gamma_{2} ; \gamma_{3} \ldots \gamma_{k}}^{\:\:\:a}\mathring{\omega}_{\bar{\alpha}_{1}\:\bar{\alpha}_{2};\bar{\alpha}_{3}\ldots\bar{\alpha}_{k}}^{\:\:\:\bar{b}} - 
g_{a\bar{b}}\omega_{\gamma_{1} \: \gamma_{2} ; \gamma_{3} \ldots \gamma_{k}}^{\:\:\:a}\omega_{\bar{\alpha}_{1}\:\bar{\alpha}_{2};\bar{\alpha}_{3}\ldots\bar{\alpha}_{k}}^{\:\:\:\bar{b}} \equiv \:\: 0.
\end{equation*}
By noticing that $(N-l') - (n - l) = (N - n) - (l' - l)$ and performing a similar argument we use Lemma 2.2 to obtain the desired equality.  The details are left to the reader.

Now we will show that the conformal equivalence is actually an equality in the expression
\begin{equation}
g_{a\bar{b}}\omega_{\gamma_{1} \: \gamma_{2} ; \gamma_{3} \ldots \gamma_{j}}^{\:\:\:a}\omega_{\bar{\alpha}_{1}\:\bar{\alpha}_{2};\bar{\alpha}_{3}\ldots\bar{\alpha}_{k}}^{\:\:\:\bar{b}} \equiv \:\: \mathring{g}_{a\bar{b}}\mathring{\omega}_{\gamma_{1} \: \gamma_{2} ; \gamma_{3} \ldots \gamma_{j}}^{\:\:\:a}\mathring{\omega}_{\bar{\alpha}_{1}\:\bar{\alpha}_{2};\bar{\alpha}_{3}\ldots\bar{\alpha}_{k}}^{\:\:\:\bar{b}}
\end{equation}
where without loss of generality, we assume $j > k$. We first assume $(N_{0} - n) + (l' - l) < l$ and rewrite equation (7.8) as
\begin{multline}
-\sum_{a=1}^{l' - l}\omega_{\gamma_{1} \: \gamma_{2} ; \gamma_{3} \ldots \gamma_{j}}^{\:\:\:a}\omega_{\bar{\alpha}_{1}\:\bar{\alpha}_{2};\bar{\alpha}_{3}\ldots\bar{\alpha}_{k}}^{\:\:\:\bar{a}}
-\sum_{b=1}^{N_{0} - n}\mathring{\omega}_{\gamma_{1} \: \gamma_{2} ; \gamma_{3} \ldots \gamma_{j}}^{\:\:\:b}\mathring{\omega}_{\bar{\alpha}_{1}\:\bar{\alpha}_{2};\bar{\alpha}_{3}\ldots\bar{\alpha}_{k}}^{\:\:\:\bar{b}} \\
+ \sum_{c=l'-l + 1}^{N-n}\omega_{\gamma_{1} \: \gamma_{2} ; \gamma_{3} \ldots \gamma_{j}}^{\:\:\:c}\omega_{\bar{\alpha}_{1}\:\bar{\alpha}_{2};\bar{\alpha}_{3}\ldots\bar{\alpha}_{k}}^{\:\:\:\bar{c}} \equiv 0.
\end{multline}
We apply a lemma of D'Angelo (see \cite{DA}, chapter 5) to equation (7.7) to obtain the existence of a unitary matrix $U$ such that
\begin{equation*}
U \left( \begin{array}{c}
\omega_{\bar{\alpha}_{1}\:\:\:\:\bar{\alpha}_{2};\bar{\alpha}_{3}\ldots\bar{\alpha}_{k}}^{\:\:\:\bar{1}} \\
\vdots \\
\omega_{\bar{\alpha}_{1}\:\:\:\:\bar{\alpha}_{2};\bar{\alpha}_{3}\ldots\bar{\alpha}_{k}}^{\:\:\:\overline{l' - l}} \\
\mathring{\omega}_{\bar{\alpha}_{1}\:\:\:\:\bar{\alpha}_{2};\bar{\alpha}_{3}\ldots\bar{\alpha}_{k}}^{\:\:\:\bar{1}} \\
\vdots \\
\mathring{\omega}_{\bar{\alpha}_{1}\:\:\:\:\bar{\alpha}_{2};\bar{\alpha}_{3}\ldots\bar{\alpha}_{k}}^{\:\:\:\overline{N_{0} - n}}\\
0\\
\vdots\\
0
\end{array} \right) = \left(
\begin{array}{c}
\omega_{\bar{\alpha}_{1}\:\:\:\:\bar{\alpha}_{2};\bar{\alpha}_{3}\ldots\bar{\alpha}_{k}}^{\:\:\:\overline{l' - l + 1}} \\
\vdots \\
\vdots \\
\vdots \\
\vdots \\
\vdots \\
\vdots \\
\omega_{\bar{\alpha}_{1}\:\:\:\:\bar{\alpha}_{2};\bar{\alpha}_{3}\ldots\bar{\alpha}_{k}}^{\:\:\:\overline{N-n}}
\end{array} \right)
\end{equation*}
Note that we are working at a fixed point here.  This implies the existence of constants $\bar{A}_{r}^{\:\:c}$ and $\bar{B}_{s}^{\:\:c}$, with $1 \le r \le l' - l$ and $1 \le s \le N_{0} - n$ such that
\begin{equation*}
\omega_{\bar{\alpha}_{1}\:\:\:\:\bar{\alpha}_{2};\bar{\alpha}_{3}\ldots\bar{\alpha}_{k}}^{\:\:\:\bar{c}} = \bar{A}_{r}^{\:\:c}\omega_{\bar{\alpha}_{1}\:\:\:\:\bar{\alpha}_{2};\bar{\alpha}_{3}\ldots\bar{\alpha}_{k}}^{\:\:\:\bar{r}} + \bar{B}_{s}^{\:\:c}\mathring{\omega}_{\bar{\alpha}_{1}\:\:\:\:\bar{\alpha}_{2};\bar{\alpha}_{3}\ldots\bar{\alpha}_{k}}^{\:\:\:\bar{s}},
\end{equation*}
where $l' - l + 1 \leq c \leq N - n$, and we are using the summation convention for the indices $r$ and $s$.

We substitute the above into equation (7.9) to obtain
\begin{multline*}
-\sum_{a=1}^{l' - l}\omega_{\gamma_{1} \: \gamma_{2} ; \gamma_{3} \ldots \gamma_{j}}^{\:\:\:a}\omega_{\bar{\alpha}_{1}\:\bar{\alpha}_{2};\bar{\alpha}_{3}\ldots\bar{\alpha}_{k}}^{\:\:\:\bar{a}}
-\sum_{b=1}^{N_{0} - n}\mathring{\omega}_{\gamma_{1} \: \gamma_{2} ; \gamma_{3} \ldots \gamma_{j}}^{\:\:\:b}\mathring{\omega}_{\bar{\alpha}_{1}\:\bar{\alpha}_{2};\bar{\alpha}_{3}\ldots\bar{\alpha}_{k}}^{\:\:\:\bar{b}} \\
+ \sum_{c=l'-l + 1}^{N-n}\omega_{\gamma_{1} \: \gamma_{2} ; \gamma_{3} \ldots \gamma_{j}}^{\:\:\:c}(\bar{A}_{r}^{\:\:c}\omega_{\bar{\alpha}_{1}\:\:\:\:\bar{\alpha}_{2};\bar{\alpha}_{3}\ldots\bar{\alpha}_{k}}^{\:\:\:\bar{r}} + \bar{B}_{s}^{\:\:c}\mathring{\omega}_{\bar{\alpha}_{1}\:\:\:\:\bar{\alpha}_{2};\bar{\alpha}_{3}\ldots\bar{\alpha}_{k}}^{\:\:\:\bar{s}}) \equiv 0.
\end{multline*}

We regroup the terms in this expression, which yields
\begin{multline*}
\sum_{r=1}^{l' - l}(\sum_{c=l'-l + 1}^{N-n}\omega_{\gamma_{1} \: \gamma_{2} ; \gamma_{3} \ldots \gamma_{j}}^{\:\:\:c}\bar{A}_{r}^{\:\:c} - \omega_{\gamma_{1} \: \gamma_{2} ; \gamma_{3} \ldots \gamma_{j}}^{\:\:\:r})\omega_{\bar{\alpha}_{1}\:\bar{\alpha}_{2};\bar{\alpha}_{3}\ldots\bar{\alpha}_{k}}^{\:\:\:\bar{r}}\\
+ \sum_{s=1}^{N_{0} - n}(\sum_{c=l'-l + 1}^{N-n}\omega_{\gamma_{1} \: \gamma_{2} ; \gamma_{3} \ldots \gamma_{j}}^{\:\:\:c}\bar{B}_{s}^{\:\:c} - \mathring{\omega}_{\gamma_{1} \: \gamma_{2} ; \gamma_{3} \ldots \gamma_{j}}^{\:\:\:s})\mathring{\omega}_{\bar{\alpha}_{1}\:\bar{\alpha}_{2};\bar{\alpha}_{3}\ldots\bar{\alpha}_{k}}^{\:\:\:\bar{s}} \equiv 0,
\end{multline*}
where we are not using the summation convention for the indices $r$ and $s$.  Since the number of terms in the sum on the left side in the preceding equation is strictly less than $n$, we use Lemma 2.1 in the same way that we used Lemma 2.2 previously to conclude that the conformal equivalence is in fact an equality.  We then recombine all terms to get the desired equality.  In the case where $(N - l') - (n - l) < l$, we apply the lemma of D'Angelo as above to obtain constants $\bar{\tilde{A}}_{c}^{r}$ and $\bar{\tilde{B}}_{c}^{s}$, such that 
\begin{equation*}
\omega_{\bar{\alpha}_{1}\:\:\:\:\bar{\alpha}_{2};\bar{\alpha}_{3}\ldots\bar{\alpha}_{k}}^{\:\:\:\bar{r}} = \bar{\tilde{A}}_{c}^{r}\omega_{\bar{\alpha}_{1}\:\:\:\:\bar{\alpha}_{2};\bar{\alpha}_{3}\ldots\bar{\alpha}_{k}}^{\:\:\:\bar{c}}, \qquad \textrm{and} \qquad \mathring{\omega}_{\bar{\alpha}_{1}\:\:\:\:\bar{\alpha}_{2};\bar{\alpha}_{3}\ldots\bar{\alpha}_{k}}^{\:\:\:\bar{s}} = \bar{\tilde{B}}_{c}^{s}\omega_{\bar{\alpha}_{1}\:\:\:\:\bar{\alpha}_{2};\bar{\alpha}_{3}\ldots\bar{\alpha}_{k}}^{\:\:\:\bar{c}},
\end{equation*}
where $1 \le r \le l' - l$, $1 \le s \le N_{0} - n$, $l' - l + 1 \leq c \leq N - n$ and we are using the summation convention on the indices $r$ and $s$.  We then substitute into (7.9) as before to obtain the desired result.  The details of this are left to the reader.

We embed the vectors representing the second fundamental form of $f_{0}$ and its derivatives into $\mathbb{C}^{N-n}$ by appending the appropriate number of zeros.  Thus we have shown that all inner products of derivatives of the second fundamental form of $f$ with respect to $g_{a\bar{b}}$ are equal to the corresponding inner products of derivatives of the second fundamental form of $f_{0}$ with respect to $\mathring{g}_{a\bar{b}}$.  Lemma 6.5 gives the desired inequality relating the dimension of $E_{k}$ and $\mathring{E}_{k}$.

In the side reversing case, the argument is similar except that we need only consider the analogue of the negative of equation (7.6).  This is because $min(N - l' - l, l' - (n-l)) = l' - (n-l)$.  We leave the details to the reader.

\renewcommand\refname{References}

Peter Ebenfelt, Ravi Shroff: Department of Mathematics, University of California at San Diego, La Jolla, CA 92093-0112, USA \\
\emph{Email address:} \textbf{pebenfel@math.ucsd.edu, rshroff@math.ucsd.edu}.

\end{document}